\documentclass[11pt]{article}
\usepackage{amsmath,amsthm,amsfonts}

\begin{document}

\newcommand{\R}{{\mathbb{R}}}

\newcommand{\Q}{{\mathbb{Q}}}
\newcommand{\T}{{\mathbb{T}}}
\newcommand{\Z}{{\mathbb{Z}}}
\newcommand{\N}{{\mathbb{N}}}
\newcommand{\C}{{\mathbb{C}}}

\newcommand{\La}{{\Lambda}}
\newcommand{\tLa}{{\widetilde{\Lambda}}}
\newcommand{\tP}{{\widetilde{\mathcal{P}}}}

\newcommand{\hJ}{{\hat{J}}}

\newcommand{\cA}{{\mathcal{A}}}
\newcommand{\cD}{{\mathcal{D}}}
\newcommand{\cE}{{\mathcal{E}}}
\newcommand{\cF}{{\mathcal{F}}}
\newcommand{\cG}{{\mathcal{G}}}
\newcommand{\cH}{{\mathcal{H}}}
\newcommand{\cI}{{\mathcal{I}}}
\newcommand{\cM}{{\mathcal{M}}}
\newcommand{\cN}{{\mathcal{N}}}
\newcommand{\cP}{{\mathcal{P}}}
\newcommand{\cS}{{\mathcal{S}}}
\newcommand{\cT}{{\mathcal{T}}}
\newcommand{\cU}{{\mathcal{U}}}
\newcommand{\cV}{{\mathcal{V}}}
\newcommand{\cW}{{\mathcal{W}}}

\newcommand{\oQ}{{\overline{QH}_{ev} (M)}}

\newcommand{\tG}{{\widetilde{G}}}
\newcommand{\tmu}{{\tilde{\mu}}}
\newcommand{\tf}{{\tilde{f}}}
\newcommand{\tlambda}{{\tilde{\lambda}}}
\newcommand{\tphi}{{\tilde{\phi}}}
\newcommand{\tih}{{\tilde{h}}}
\newcommand{\g}{{\mathfrak{g}}}

\newcommand{\Ca}{{\hbox{\it Cal}}}
\newcommand{\tCa}{{\widetilde{\hbox{\it Cal}}}}

\newcommand{\Ham}{{\hbox{\it Ham}}}
\newcommand{\Symp}{{\hbox{\it Symp}}}
\newcommand{\spec}{{\text{spec}}}
\newcommand{\id}{{\text{{\bf 1}}}}
\def\square{{\vrule height6pt width6pt depth2pt}}

\newtheorem{thm0}{Theorem}[section]
\newtheorem{thm}{Theorem}[subsection]
\newtheorem{thm1}{Theorem}[subsubsection]
\newtheorem{defin0}[thm0]{Definition}
\newtheorem{defin}[thm]{Definition}
\newtheorem{defin1}[thm1]{Definition}
\newtheorem{cor}[thm]{Corollary}
\newtheorem{cor1}[thm1]{Corollary}
\newtheorem{lem}[thm]{Lemma}
\newtheorem{lem1}[thm1]{Lemma}
\newtheorem{prop}[thm]{Proposition}
\newtheorem{rem}[thm]{Remark}
\newtheorem{exam}[thm]{Example}
\newtheorem{exam1}[thm1]{Example}
\newtheorem{sublem}[thm]{Sublemma}
\newtheorem{assum}[thm]{Assumption}
\def\b{{\vrule height6pt width6pt depth2pt}}

\newcommand{\tHam}{\widetilde{\hbox{\it Ham}\, }}
\newcommand{\tSymp}{\widetilde{\Symp}}

\title{Commutator length of symplectomorphisms}

\author{
Michael Entov\thanks{Partially supported by the Nathan, Guta and
Robert Gleit\-man Post\-doc\-toral Fellow\-ship and by the Israel
Science Foundation founded by the Israel Academy of Sciences and
Humanities, grant n. 582/99-1}
\\
\\
Department of Mathematics\\ Technion -- Israel Institute of
Technology\\ Haifa 32000, Israel\\
e-mail: entov@math.technion.ac.il\\}

\date{\today}
\maketitle

\begin{abstract}

Each element $x$ of the commutator subgroup $[G, G]$ of a group $G$
can be represented as a product of commutators. The minimal number
of factors in such a product is called the commutator length of
$x$. The commutator length of $G$ is defined as the supremum of
commutator lengths of elements of $[G, G]$.

We show that for certain closed symplectic manifolds $(M,\omega)$,
including complex projective spaces and Grassmannians,
the universal cover $\tHam  (M,\omega)$ of the group of Hamiltonian
symplectomorphisms of  $(M,\omega)$ has infinite
commutator length.  In particular, we present explicit examples of
elements in $\tHam  (M,\omega)$ that have
arbitrarily large
commutator length -- the estimate on their commutator length
depends on the multiplicative structure of the quantum cohomology
of $(M,\omega)$. By a different method we also show that in the
case $c_1 (M) = 0$ the group
$\tHam  (M,\omega)$ and the universal cover
$\tSymp_0\, (M,\omega)$ of the identity
component of the group of symplectomorphisms of $(M,\omega)$
have infinite commutator length.

\end{abstract}



\vfil \eject

\tableofcontents

%\bigskip
%\bigskip
\section{An overview of the main results}
\label{sec-commut-lengths-defs-main-results}

\subsection{Basic definitions and preliminaries}

Let $G$ be a group. Each element $x$ of the commutator subgroup
$[G,G]\subset G$ can be represented as a product of commutators.
The minimal number of factors in such a product is called the {\it
commutator length of $x$} and will be denoted by ${\hbox{\it cl}}\,(x)$
(sometimes this quantity is also called "the genus of $x$").
Define the {\it commutator length of $G$} as ${\hbox{\it cl}}\,
(G) := \sup_{x\in [G, G]}\, {\hbox{\it cl}}\,(x)$. Define the {\it
stable commutator norm} $\| x\|_{cl}$ of $x\in [G,G]$ as $\|
x\|_{cl} := \lim_{k\to +\infty} {\hbox{\it cl}}\, (x^k)/k$. A
simple exercise shows that such a
limit always exists.

For results on commutator length of finite-dimensional Lie groups
see e.g. \cite{Bor}, \cite{Goto}. In this paper we will study the
commutator lengths of the universal covers of the
infinite-dimensional Lie groups $\Symp_0\, (M,\omega)$
and $\Ham\, (M,\omega)$, associated with a closed
connected symplectic manifold $(M, \omega)$.

Here are the definitions of these groups. Let $(M,\omega)$
be a closed connected symplectic manifold.
By $\Symp_0\, (M,\omega)$ we denote the identi\-ty
component of the group of symplectomorphisms of $(M,\omega)$. Its
universal cover will be denoted by
$\tSymp_0\, (M,\omega)$.

The group $\Ham\, (M,\omega)$ of {\it Hamiltonian
symplectomorphisms} of $(M,\omega)$ is defined as follows. A
function $H: S^1\times M\to \R$ (called {\it Hamiltonian})
defines a time-dependent {\it Hamiltonian vector
field} $X (t,x)$, $(t,x)\in S^1\times M$, on $M$ by the
formula\footnote{Note our sign convention. Also note the
signs in the formula (\ref{eqn-def-action})
for the action functional and in the formula (\ref{eqn-cz-index})
for the Conley-Zehnder index.}
\begin{equation}
\label{eqn-def-ham-vf} d H^t (\cdot) = \omega (X^t, \cdot),
\end{equation}
where $H^t  = H (t,\cdot)$, $X^t = X (t, \cdot)$, and the
formula holds pointwise on $M$ for every $t\in S^1$.
 The flow of a Hamiltonian vector
field preserves the symplectic form. A {\it Hamiltonian
symplectomorphism} $\varphi$ of $(M,\omega)$ is a symplectomorphism
of $M$ which can be represented as the time-1 map of the flow
(called {\it Hamiltonian flow}) of a (time-dependent) Hamiltonian
vector field. If this Hamiltonian vector field is defined by a
Hamiltonian function $H$ we say that $\varphi$ and the whole
Hamiltonian flow are {\it generated} by $H$.

Observe that different Hamiltonian functions can generate the same
(time-dependent) Hamiltonian vector field and therefore the same
Hamiltonian flow and the same Hamiltonian symplectomorphism.
Namely, given a Hamiltonian vector field $X$ generated by a
Hamiltonian $H: S^1\times M\to \R$ one can write any other
Hamiltonian generating $X$ as $H^\prime = H + h (t)$ for some
function $h: S^1\to \R$. In order to deal with this
non-uniqueness we introduce the following normalization.

%\bigskip
%\noindent
\begin{defin}
\label{def-normalized-Hamilt-n} {\rm Suppose $M$ is closed. A
Hamiltonian $H:S^1\times M\to \R$ is called {\it normalized}
if $\int_M H^t \omega^n = 0$, for any $t\in S^1$. }
\end{defin}
%\bigskip

Thus any Hamiltonian vector field (or Hamiltonian flow) is
generated by a unique normalized Hamiltonian function.

Hamiltonian symplectomorphisms of $(M,\omega)$ form a group
$\Ham\, (M,\omega)\subset \Symp_0\,
(M,\omega)$. In particular, if Hamiltonian flows $\{ \varphi^t\},
\{\psi^t\}$ are generated, respectively, by Hamiltonians $F ,
G$, then the flow $\{ \varphi^t \cdot\psi^t\}$ is generated by
the Hamiltonian function $F\sharp G (t,x) := F (t, x) + G (t,
(\varphi^t)^{-1} (x))$ and the flow $\{ (\varphi^t)^{-1}\}$ is
generated by the Hamiltonian function $\bar{F} (t,x) = - F
(t,\varphi^t (x))$.

Denote by $\tHam  (M,\omega)$ the universal cover of $\Ham\,  (M,\omega)$.
One can view $\tHam  (M,\omega)$ as a
subgroup of $\tSymp_0\, (M,\omega)$ \cite{Ban},
\cite{Cal}. A theorem of A.Banyaga \cite{Ban} shows that
\begin{eqnarray}
\label{eqn-ham-symp-2} \tHam  (M,\omega) &
= & [\tHam  (M,\omega),
\tHam  (M,\omega)] = \nonumber \\ & &
[\tSymp_0\, (M,\omega), \tSymp_0\, (M,\omega)].
\end{eqnarray}

%\smallskip
\begin{defin}
\label{def-elems-of-ham-tilde-ham-generated-by-H} {\rm We will
denote by $\varphi_H$ the Hamiltonian symplecto\-mor\-phism generated by
a Hamiltonian function $H$. The Hamiltonian flow generated by $H$ is
a path connecting $Id$ and $\varphi_H$ which defines a lift of
$\varphi_H$ in $\tHam  (M,\omega)$. This lift
will be denoted by $\widetilde{\varphi}_H$.
}
\end{defin}
%\bigskip

Obviously, $\widetilde{\varphi}_F^{-1} =
\widetilde{\varphi}_{\bar{F}}$, $\widetilde{\varphi}_{F} \cdot
\widetilde{\varphi}_{G} = \widetilde{\varphi}_{F\sharp G}$
and if $F$ is time-independent then $\widetilde{\varphi}_{kF} =
\widetilde{\varphi}_{F}^k$ for any $k\in\Z$.

A {\it ball} in $(M^{2n}, \omega)$ is the image of a symplectic
embedding into $M$ of a round $2n{\hbox{\rm -dimensional}}$
open ball in the standard
symplectic $\R^{2n}$.
We say that a
ball $B\subset M$ is {\it displaced} by
$\varphi\in \Ham\,  (M,\omega)$ if
\[
\varphi (B) \cap {\hbox{\rm Closure}}\, (B) = \emptyset.
\]
A ball $B\subset M$ is called {\it
displaceable} if it can be displaced by some
$\varphi\in \Ham\,  (M,\omega)$.

%\smallskip
\begin{defin}
\label{def-simple-Hamiltonians} {\rm
Let $B\subset M$ be a ball.
Identify $B$ symplectically with the ball $B (0, R)$
of a radius $R$ with the center at zero in the standard symplectic
$\R^{2n}$. Define a time-independent Hamiltonian $H_B: M\to\R$ as
follows. Outside $B$ set $H_B (x) \equiv -1$. Inside $B$ set
$H_B: B (0,R)\to \R$
as
$H_B (x) := \zeta (\| x\|)$, where $\| \cdot\|$ is the Euclidean
norm on $\R^{2n}$ and $\zeta: \R\to \R$
is a smooth function which is equal to a positive constant
near $0$ and to $-1$ outside of $[-R/2,
R/2]$, and which is chosen in such a way that
$\int_M H_B\, \omega^n = 0$ (thus $H$ is a normalized
Hamiltonian). }
\end{defin}
%\bigskip

Finally, we say that an almost complex structure $J$ on
$(M,\omega)$ is {\it compatible} with $\omega$ if $\omega (\cdot,
J\cdot)$ is a Riemannian metric on $M$. Almost complex structures
compatible with $\omega$ form a non-empty contractible space
\cite{Gro-pshc}. The first Chern class of $T M$ for any $J$
compatible with $\omega$ is the same and will be denoted by $c_1
(M)$.

%\bigskip
\subsection{The results}
\label{subsect-results-nice-form}

Let $M = \C P^n$ be the complex projective space equipped with
the standard Fubini-Study symplectic form $\omega$.
Let $A = [\C P^1]\in H_2 (\C P^n)$ be the homology class of a
projective line in $\C P^n$.

%\smallskip
\begin{thm}
\label{thm-cl-cp-n}
Let $B\subset \C P^n$ be a ball and let $F: S^1\times \C P^n\to\R$
be a normalized
Hamiltonian so that $\varphi_F$ displaces $B$.
Let $g$ be a positive
integer. Then for any
\[
w > \int_0^1 \sup_M F (t, \cdot  ) dt + \Bigg[ \frac {gn}{n+1}
\Bigg] \omega (A)
\]
one has ${\hbox{\it cl}}\,(\widetilde{\varphi}_{F\sharp w H_B}) >
g$ in the group $\tHam  (\C P^n,\omega)$. Hence
\[
{\hbox{\it cl}}\, (\tHam  (\C P^n,\omega)) = +\infty.
\]
\end{thm}
%\smallskip

The statement ${\hbox{\it cl}}\, (\tHam  (\C P^n,\omega)) = +\infty$ can be strengthened: one can show that the
infinite supremum of commutator lengths is actually reached on
some cyclic subgroup of $\tHam  (\C P^n,\omega)$.

%\smallskip
\begin{thm}[L.Polterovich]
\label{thm-stable-norm-cp-n} Let
$B\subset \C P^n$ be a displaceable ball. Then for any $w > 0$
\[
\| \widetilde{\varphi}_{wH_B} \|_{cl} \geq \frac{w (n+1)}{n \omega
(A)} > 0.
\]
\end{thm}
%\smallskip

A closed connected symplectic manifold $(M, \omega)$ is
called {\it symplectically aspherical} if $\omega$ vanishes on any
spherical homology class.

%\smallskip
\begin{thm}[M.Entov-L.Polterovich]
\label{thm-ostrover-examples-sympl-aspher} Let $(M, \omega)$ be
symplectically aspherical. Let $B\subset \C P^n$ be a ball
and let $F: S^1\times \C P^n\to\R$
be a normalized
Hamiltonian so that $\varphi_F$ displaces $B$.
Then for any
\[
w > \int_0^1 \sup_M F (t, \cdot  ) dt
\]
one has ${\hbox{\it cl}}\,(\widetilde{\varphi}_{F\sharp w H_B}) >
1$ in the group $\tHam  (M,\omega)$.

\end{thm}
%\smallskip

A closed connected symplectic surface $(M,\omega)$ of positive
genus is  symplectically aspherical and $\Ham\,  (M,\omega) = \tHam  (M,\omega)$
(see \cite{Pol-book}). In this case
Theorem~\ref{thm-ostrover-examples-sympl-aspher} provides
examples of elements of $\Ham\,  (M,\omega)$ with
commutator length greater than 1.

Theorem~\ref{thm-cl-cp-n} has a direct analogue for complex
Grassmannians. We will state here only the following conclusion.
Let $M = Gr\, (r, n)$, $1\leq r\leq n-1$, be the Grassmannian of
complex $r{\hbox{\rm -dimensional}}$ subspaces in $\C^n$. It
carries a natural symplectic form
$\omega$ (see Example~\ref{exam-grassm}). Let $A$ be the
generator of $H_2 (Gr\, (r, n), \Z)\cong \Z$ such that $\omega (A) >0$.

%\smallskip
\begin{thm}
\label{thm-cp-n-grassms-infinite-cl} For a displaceable ball
$B\subset Gr\, (r, n)$  and any $w > 0$
\[
\| \widetilde{\varphi}_{wH_B} \|_{cl} \geq \frac{w r (n-r)}{n
\omega (A)} > 0.
\]
in $\tHam  (Gr\, (r, n),\omega)$. Hence
${\hbox{\it cl}}\,(\tHam  (Gr\, (r,
n),\omega)) = +\infty$.
\end{thm}
%\smallskip

The proofs of the results above can be found in
Section~\ref{sect-pfs-thm-cl-cp-n-thm-stable-norm-cp-n}.
They are based on the Floer theory.

Observe that these results provide an
estimate on the commutator length of a certain {\it individual}
element $\varphi\in \tHam  (M,\omega)$ and the estimates on
$\|\varphi\|_{cl}$ are obtained as a {\it consequence of estimates
on the commutator length of each individual} $\varphi^k$, $k\in\Z$.
Another way of estimating $\|\varphi\|_{cl}$ is related to the
notion of {\it a quasimorphism}. This approach gives no
information on the number ${\hbox{\it cl}}\, (\varphi^k)$ for an {\it
individual} $k$ -- it reflects only the {\it asymptotics} of those
numbers as $k\to +\infty$.

Recall that {\it a quasimorphism} on a group $G$ is a function
${\mathfrak f}: G \to \R$ which satisfies the homomorphism
equation up to a bounded error: there exists $R > 0$ such that
$$|{\mathfrak f}(ab) -{\mathfrak f}(a) -{\mathfrak f}(b)| \leq R$$
for all $a,b \in G$. A quasimorphism ${\mathfrak f}$ is called
{\it homogeneous} if ${\mathfrak f} (a^m) = m {\mathfrak f} (a)$
for all $a \in G$ and $m \in \Z$. A simple  exercise shows that if
there exists a homogeneous quasimorphism which does not vanish on
$\phi\in [G,G]$ then $\| \phi\|_{cl} >0$ and hence ${\hbox{\it
cl}}\,(G) = +\infty$. Conversely, according to a theorem of
C.Bavard \cite{Bav}, if $\| \phi\|_{cl} >0$ then there exists a
homogeneous quasimorphism on $G$ which does not vanish on $\phi$.

%\smallskip
\begin{thm}
\label{thm-comm-length-closed-mfds} Let $(M,\omega)$ be a closed
connected symplectic manifold with $c_1 (M) = 0$ and let $G =
\tSymp_0\, (M,\omega)$. Then there exists a
homogeneous quasimorphism ${\mathfrak f}: G\to \R$ which does not
vanish on $[G,G] = \tHam  (M,\omega)$ and
therefore
\[
{\hbox{\it cl}}\,(\tSymp_0\, (M,\omega)) =
{\hbox{\it cl}}\, (\tHam  (M,\omega)) =
+\infty.
\]
In particular, for any (not necessarily displaceable)
ball $B\subset M$ one has ${\mathfrak f}
(\widetilde{\varphi}_{H_B}) \neq 0$ and hence
\[
\|\widetilde{\varphi}_{H_B}\|_{cl} > 0
\]
both in $\tHam  (M,\omega)$ and in $\tSymp_0 (M,\omega)$.
\end{thm}
%\smallskip

%\smallskip
\begin{exam}
\label{exam-c-1-zero} {\rm Here are examples of closed
symplectic manifolds with $c_1 (M) = 0$:

%\medskip
\noindent 1) Symplectic tori.

%\smallskip
\noindent 2) Complex K{\"a}hler manifolds with $c_1 = 0$. These
are Ricci-flat mani\-folds that include complex K3-surfaces,
hyper-K{\"a}hler manifolds, Calabi-Yau manifolds of complex
dimension 3 and certain smooth complete intersections in $\C P^n$, $n\geq 2$. The last example can be explicitly described as
follows. If $M\subset \C P^n$ is a smooth complex submanifold
of complex dimension $k$, which is the intersection of $n-k$
hypersurfaces of degrees $d_1,\ldots, d_{n-k}$, then $c_1 (M) = 0$
as soon as $\sum_{i=1}^{n-k} d_i = n+1$.

}
\end{exam}
%\smallskip

Theorem~\ref{thm-comm-length-closed-mfds} is proven in
Section~\ref{sect-quasimorphisms}. Our construction of ${\mathfrak
f}$ directly generalizes the one of J.Barge and E.Ghys
\cite{Ba-Ghys} who proved a similar result for the group
$\Symp^{c}\, (B^{2n})$ of compactly supported
symplectomorphisms of the standard symplectic ball $B^{2n}$.    A
similar construction (albeit used for different purposes)
is presented in
\cite{Co-Ga-It-Pa}.

%\bigskip
\begin{rem}
\label{rem-other-quasimorphisms} {\rm

1) According to the theorem of C.Bavard \cite{Bav} mentioned
above, Theorems~\ref{thm-stable-norm-cp-n} and
\ref{thm-cp-n-grassms-infinite-cl} yield the existence of a
homogeneous quasimorphism on $\tHam  (M,\omega)$ for complex
projective spaces and Grassmannians. However the proof of
Bavard's result \cite{Bav} uses Hahn-Banach theorem and does not
provide any explicit construction of a quasimorphism on
$\tHam  (M,\omega)$. Such an explicit construction has been carried
out in \cite{EP} for complex
projective spaces, Grassmannians and certain other symplectic
manifolds. This in turn yields an alternative proof of
Theorems~\ref{thm-stable-norm-cp-n} and
\ref{thm-cp-n-grassms-infinite-cl}. In fact for complex projective
spaces the quasimorphism on $\tHam  (\C P^n,\omega)$ constructed
in \cite{EP}
descends to a homogeneous quasimorphism on $\Ham\,  (\C P^n,\omega)$,
yielding, in particular,
\[
\| \varphi_{H_B}\|_{cl} >0
\]
in $\Ham\,  (\C P^n,\omega)$ for a displaceable ball
$B\subset \C P^n$ (see \cite{EP}). Hence
\[
{\hbox{\it cl}}\, (\Ham\,  (\C P^n,\omega)) = +\infty.
\]

\bigskip
\noindent 2) Using completely different methods J.-M.Gambaudo and
E.Ghys \cite{GG} constructed homogeneous quasimorphisms on
$G=\Symp^{c}_0\, (M,\omega)$ which do not vanish on
$[G,G] = \Ham^{c}\, (M,\omega)$ for any two-dimensional
symplectic manifold $(M,\omega)$. Thus the commutator lengths of
both groups $\Symp^{c}_0\, (M,\omega)$, $\Ham^c\,
(M,\omega)$ are infinite. }
\end{rem}

%\bigskip
\noindent {\bf Acknowledgments. } The original idea of applying
K-area to the study of commutator length, which inspired the key
result in Section~\ref{subsect-def-K-area} of this paper, belongs to
F.Lalonde. I am most grateful to L.Polterovich for suggesting the
problem to me, for important comments, corrections and
numerous extremely useful discussions. I thank D.Fuchs for
explaining to me the simple proof of Lemma~\ref{lem-fuchs}
presented in this paper. I also thank A.Postnikov for
communicating to me Proposition~\ref{prop-post}. I am obliged to
G.Lu, Y.-G. Oh and the anonymous referee who read a preliminary
version of this paper and made valuable comments and corrections.
This work was started during my stay at the Weizmann Institute of
Science -- I thank it for the warm hospitality.

%\bigskip
%\bigskip
\section{Commutator length of elements in
$\tHam  (M,\omega)$ }
\label{sec-explicit-symps-with-cl-geq-2}

In this chapter we state general results concerning estimates on
the commutator length of individual elements of
$\tHam  (M,\omega)$ in terms of the Floer
theory and give examples of such elements with a large
commutator length. In
Sections~\ref{subsec-str-semipos-symp-mfds}-\ref{subsect-spec-numbers-intro}
we introduce the necessary preliminaries about quantum and Floer
(co)homology. The main results are stated in
Sec\-tions~\ref{subsec-thm-action-length}. In
Section~\ref{subsec-thm-action-length-applications} we state practical
applications of the main results and use them to prove
Theorems~\ref{thm-cl-cp-n}-\ref{thm-cp-n-grassms-infinite-cl}.

%\bigskip
\subsection{Strongly semi-positive symplectic manifolds}
\label{subsec-str-semipos-symp-mfds}

Let $\Pi$ denote the group of non-torsion spherical homology
classes of $M$, i.e. the image of the Hurewicz homomorphism $\pi_2
(M)\to H_2 (M, {\bf Z})/Tors$.

A closed symplectic manifold $(M^{2n},\omega)$
is called
{\it strongly semi-positive} \cite{ME}, if for every
$A\in \Pi$
\[
2-n\leq c_1 (A)< 0 \Longrightarrow \omega (A)\leq 0.
\]
In \cite{Se} such manifolds were said to satisfy the "Assumption
$W^{+}$".

As in \cite{ME}, for technical reasons
the closed symplectic manifold $(M,\omega)$ is assumed to be
strongly semi-positive whenever we use Floer or quantum (co)homology,
though in view of the recent
developments (see \cite{Fu-Ono}, \cite{Liu-Ti}, \cite{Liu-Ti1},
\cite{Lu}) it is likely that this assumption can be removed.

The class of strongly semi-positive manifolds includes all
symplectic manifolds of dimension less or equal to 4.
A particularly important subclass of strongly semi-positive manifolds
is formed by {\it spherically monotone} symplectic manifolds, i.e. the ones
with $\left. c_1\right|_\Pi = \kappa \left. [\omega]\right|_\Pi$
for some $\kappa >0$.
Complex projective spaces,
complex Grassmannians and flag spaces and symplectically aspherical
manifolds are spherically monotone.

%\bigskip
\subsection{Novikov ring, quantum cohomology and the Euler class}
\label{subsec-Novikov-quantum}

The {\it rational Novikov ring} $\Lambda_\omega$ of a symplectic
manifold
$(M,\omega)$ is built from the rational group ring of the
group $\Pi$ by means of the symplectic form $\omega$ \cite{Ho-Sa}.
Namely, an element of $\Lambda_\omega$ is a formal sum
\[
\lambda = \sum_{A\in\Pi} \lambda_A e^{A}
\]
with rational coefficients $\lambda_A\in \Q$ which satisfies
the condition
\[
\sharp \{ A\in\Pi\, |\, \lambda_A\neq 0, \omega (A) \leq c\}
<\infty \ \ {\rm for\ every}\ c > 0.
\]
The natural multiplication turns
$\Lambda_\omega$ into a ring.

The grading on $\Lambda_\omega$,
defined by the
condition ${\rm deg}\, (e^{A}) = 2 c_1 (A)$, fits with the ring
structure.
Sometimes we will need to consider the opposite grading on
$\Lambda_\omega$ so that the degree of $(e^A)$ is $- 2 c_1 (A)$.
The latter
grading does not fit with the multiplicative structure of
$\Lambda_\omega$. When using the second grading we will say
that $\Lambda_\omega$ is {\it anti-graded}, while using
the first grading we will call $\Lambda_\omega$ simply {\it
graded}.
The set of elements of
degree $i$ of the {\it graded} $\Lambda_\omega$ will be denoted by
$\Lambda^i_\omega$.

The {\it quantum cohomology} of $M$ is a
module over $\Lambda_\omega$ defined as a graded tensor product
$QH^\ast (M,\omega) = H^\ast (M, \Q) \otimes_{\Q} \Lambda_\omega$
over $\Q$ (here $\Lambda_\omega$ is taken in the graded version).
More importantly the group $QH^\ast (M,\omega)$ carries a delicate
multiplicative structure
called the {\it quantum product} -- we refer the reader to
\cite{Ru-Ti}, \cite{Ru-Ti-1}, \cite{Wi} for the
definitions. The quantum product is associative,
$\Lambda_\omega{\hbox{\rm -linear}}$ and respects the grading on
$QH^\ast (M,\omega)$. The cohomology class ${\bf
1}\in H^0 (M)$, Poincar{\'e}-dual to the fundamental class $[M]$, is
a unit element in $QH^\ast (M,\omega)$. We will denote by $m\in
H^{2n} (M^{2n}, \Q)$ the singular cohomology class
Poincar{\'e}-dual to a point.

Similarly, {\it quantum homology} is defined as a graded tensor product
\[
QH_\ast (M,\omega) = H_\ast (M, \Q)\otimes_{\Q}\Lambda_\omega\]
over $\Q$,
where $\Lambda_\omega$ is taken in the {\it anti-graded} version.
Thus $QH_\ast (M,\omega)$ is a graded module over the anti-graded
ring $\Lambda_\omega$. There is a natural evaluation pairing
\[
(\cdot ,\cdot ): QH^k (M,\omega) \times QH_k
(M,\omega)\to\Lambda_\omega^0.\] Namely, let $\alpha e^A\in QH^k (M,
\omega)$, $\alpha\in H^i (M, \Q)$, $e^A\in \Lambda_\omega^{k-i}$,
and $\beta e^B\in QH_k (M, \omega)$, $\beta\in H_j (M, \Q)$,
$e^B\in \Lambda_\omega^{j-k}$ (recall that elements of
$\Lambda_\omega^{j-k}$ have degree $k-j$ in the anti-graded
version of $\Lambda_\omega$). Then
\[
(\alpha e^A, \beta e^B ) = \alpha (\beta) e^{A+B},
\]
where $\alpha (\beta) \in \Q$ is the result of classical
evaluation of cohomology class on homology (in particular,
$(\alpha e^A, \beta e^B ) = 0$ if $i\neq j$ above). This is a
non-degenerate pairing. Since $\Q\subset \Lambda_\omega^0$
the universal coefficient theorem implies that the pairing $(\cdot
,\cdot )$ defines an {\it isomorphism of groups}:
\[
QH^k (M,\omega) \cong {\hbox{\rm Hom}} (QH_k (M, \omega),
\Lambda_\omega^0 ),
\]
where ${\hbox{\rm Hom}} (QH_k (M, \omega), \Lambda_\omega^0 )$ is
the group of all {\it group homomorphisms} \break $QH_k (M,
\omega)\to \Lambda_\omega^0$.

The Poincar{\'e} isomorphism $PD: QH^\ast (M,\omega) \to QH_{2n-\ast}
(M, \omega)$ is defined by the formula:
\[
PD (\sum_{A\in \Pi} \alpha_A e^A ) := \sum_{A\in \Pi} PD(\alpha_A)
e^A,\ \ \alpha_A\in H_\ast (M, \Q),
\]
where in the right-hand side $PD$ denotes the classical
Poincar{\'e} isomorphism between singular homology and cohomology.
One can check that $PD: QH^\ast (M,\omega) \to QH_{2n-\ast} (M,
\omega)$ is the homomorphism of $\Lambda_\omega{\hbox{\rm
-modules}}$ (recall that $QH_\ast (M, \omega)$ is a graded module
over the {\it anti-graded} ring $\Lambda_\omega$).

The Poincar{\'e} isomorphism and the evaluation  pairing $(\cdot,
\cdot)$ together define a nondegenerate pairing
$\langle\cdot ,\cdot\rangle$:
\[
\langle \cdot , \cdot \rangle : QH^\ast (M,\omega)\times QH^{2n-\ast}
(M,\omega)\to \Lambda_\omega^0,
\]
where
\[
\langle \alpha , \beta\rangle = ( \alpha, PD (\beta) ).
\]
In particular, if $\alpha\in QH^i (M,\omega)$, $\beta\in QH^j
(M,\omega)$ and $i+j\neq 2n$ then $\langle \alpha , \beta\rangle =
0$.

The pairing $\langle \cdot , \cdot \rangle$ can be also expressed
in terms of the quantum multiplication (see \cite{PSS},
\cite{Sch-PhD}, cf.
Section~\ref{sect-pss-stuff}).
Namely, for any $\alpha,\beta \in QH^\ast (M,\omega)$
\[
\langle\alpha, \beta \rangle = (\alpha\ast \beta, [M]),
\]
where $\alpha\ast \beta \in QH^\ast (M, \omega)$ is the quantum
product of $\alpha$ and $\beta$ and the fundamental class $[M]$ is
viewed as an element of $QH_{2n} (M, \omega)$. In other words, if
$\alpha\ast \beta$ is represented as
$\alpha\ast \beta= \sum_b \lambda_b b$, $b\in H^\ast
(M, \Q)$, $\lambda_b\in \Lambda_\omega$, then $\langle\alpha, \beta \rangle$
is the $\Lambda_\omega^0{\hbox{\rm -component}}$ of the coefficient
$\lambda_m\in\Lambda_\omega$ at the singular cohomology class $m$.

Pick a basis $\{ e_i\}$, $e_i\in H^{k_i} (M, \Q)$, of $H^\ast
(M, \Q)$ over $\Q$. Then $\{ e_i\}$
is also a basis of
$QH^\ast (M,\omega)$ over $\Lambda_\omega$. Consider the basis $\{
F_i \}$, $F_i\in H_i (M, \Q)$,
of $H_\ast (M, \Q)$ over $\Q$ dual to $\{ e_i\}$, i.e. $e_i (F_j) =
\delta_{ij}$. Then $\{ F_i \}$ is also a basis of $QH_\ast
(M,\omega)$ over $\Lambda_\omega$. Set $\bar{e}_i = PD (F_i)\in
QH^{2n-i} (M, \omega)$ for every $i$.
The isomorphism $PD$ is $\Lambda_\omega{\hbox{\rm
-linear}}$ and therefore $\{ \bar{e}_i\}$ is a basis of $QH^\ast
(M,\omega)$ over $\Lambda_\omega$. Thus
\begin{equation}
\label{eqn-poincare-pairing-of-bases} \langle e_i, \bar{e}_j
\rangle = \delta_{ij}.
\end{equation}

Any two $\Lambda_\omega{\hbox{\rm -bases}}$ $\{ e_i\}$,
$\{ \bar{e}_i\}$,
$e_i\in QH^{k_i} (M,\omega)$, $\bar{e}_i\in QH^{2n-k_i} (M,\omega)$,
of $QH^\ast (M,\omega)$ satisfying
(\ref{eqn-poincare-pairing-of-bases}) will be called
{\it Poincar{\'e}-dual } to each other.
Given such bases $\{ e_i\}$, $\{ \bar{e}_i\}$, set
\[
E=\sum_i (-1)^{k_i} e_i\ast \bar{e}_i.
\]
One can easily check
that $E\in QH^{2n} (M, \omega)$, called the {\it Euler class},
does not depend on the choice of the Poincar{\'e}-dual
bases $\{ e_i\}$, $\{
\bar{e}_i\}$.
The component of degree $2n$ of $E = \sum_i (-1)^{{\hbox{\rm
deg}}\, e_i}  e_i\ast \bar{e}_i$ is equal to $\chi (M) m\in
H^{2n} (M, \Q)$, where $\chi (M)$ is the Euler characteristic.

%\bigskip
\subsection{Floer cohomology: basic definitions}
\label{subsect-floer-defs}

Consider pairs $(\gamma, f)$, where $\gamma: S^1\to M$ is a
contractible curve and $f: D^2\to M$ is a disk spanning
$\gamma$, i.e. ${\left. f\right|}_{\partial D^2} = \gamma$.
Two pairs $(\gamma, f)$ and $(\gamma, f^\prime )$ are
called {\it equivalent } if the connected sum $f\sharp f^\prime$
represents a torsion class in $H_2 (M, {\bf Z})$.

Given a (time-dependent) Hamiltonian function $H$ on $M$ denote by
${\cal P} (H)$ the space of equivalence classes ${\hat\gamma} =
[\gamma, f ]$ of pairs $(\gamma, f)$, where $\gamma: S^1\to M$ is
a contractible time-1 periodic trajectory of the Hamiltonian flow
of $H$. The theorems proving Arnold's conjecture imply that
any Hamiltonian symplectomorphism of $M$ has a closed contractible
1-periodic orbit and therefore ${\cal P} (H)$ is always non-empty
(see \cite{Fu-Ono}, \cite{Liu-Ti}; in the semi-positive case the
conjecture was first proven in \cite{Ho-Sa}).

For an element $\hat{\gamma} = [\gamma, f]\in {\cal P} (H)$ define
its {\it action: }
\begin{equation}
\label{eqn-def-action} {\cal A}_H ({\hat\gamma}) := - \int_D
f^\ast \omega - \int_0^1 H (t, \gamma (t)) dt.
\end{equation}

If one views ${\cal A}_H$ as a functional on the space of all
equivalence classes $[\gamma, f ]$ for all contractible loops
$\gamma$ then ${\cal P} (H)$ is
precisely the set of critical points of ${\cal A}_H$ on such a space.

The group $\Pi$ acts
on ${\cal P} (H)$: if $\hat{\gamma} = [\gamma, f]\in {\cal P}
(H)$, $A\in\Pi$ then
\begin{equation}
\label{eqn-cal-P-action-of-Pi-one-Ham}
A: \hat{\gamma}\mapsto A\sharp
\hat{\gamma} := [\gamma, A\sharp f].
\end{equation}
Note that
\begin{equation}
\label{eqn-action-funct-action-of-Pi} {\cal A}_H (A\sharp
{\hat\gamma}) = {\cal A}_H ({\hat\gamma}) -\omega (A)
\end{equation}
for any ${\hat\gamma}\in {\cal P} (H)$, $A\in\Pi$.

Denote by ${\hbox{\it Spec}}\, (H)$ the set of values of ${\cal
A}_H$ on the elements of ${\cal P} (H)$.

%\smallskip
\begin{prop}[\cite{HZ}, \cite{Oh-act}, \cite{Sch-1}]
\label{prop-spectr-measure-zero} ${\hbox{\it Spec}}\, (H)\subset
\R$ is a set of measure zero.
\end{prop}

Later, in Section~\ref{sec-upsilon-widetilde-Ham-psh-curves}, we will also need the following
product version of the previous definitions. Given $l$
Hamiltonians $H = (H_1,\ldots ,H_l)$ on $M$ we denote by ${\cal P}
(H)$ the set of equivalence classes
$[{\hat{\gamma}}_1, \ldots ,{\hat{\gamma}}_l]$
of tuples
$({\hat{\gamma}}_1, \ldots ,{\hat{\gamma}}_l)$, where
${\hat{\gamma}}_i = [\gamma_i, f_i] \in {\cal P}(H_i)$ and the
equivalence relation is given by
\[
({\hat{\gamma}}_1, \ldots ,{\hat{\gamma}}_l )\sim (A_1 \sharp
{\hat{\gamma}}_1, \ldots , A_l \sharp {\hat{\gamma}}_l ),
\]
whenever $A_i\in \Pi$ and $A_1+\ldots +A_l$ is a torsion class.
The group $\Pi$ acts on ${\cal P} (H)$:
\begin{equation}
\label{eqn-cal-P-action-of-Pi} A: [{\hat{\gamma}}_1, \ldots
,{\hat{\gamma}}_l ]\mapsto [A\sharp {\hat{\gamma}}_1,
{\hat{\gamma}}_2,\ldots ,{\hat{\gamma}}_l ].
\end{equation}
For $\hat{\gamma} := [{\hat{\gamma}}_1, \ldots ,{\hat{\gamma}}_l]\in {\cal P} (H)$
set ${\cal A}_H ({\hat\gamma}) = {\cal A}_{H_1}
(\hat{\gamma}_1)+\ldots + {\cal A}_{H_l} (\hat{\gamma}_l)$.
As before,
\[
{\cal A}_H (A\sharp
\hat{\gamma}) = {\cal A}_H (\hat{\gamma}) -\omega (A)
\]
for any ${\hat\gamma}\in {\cal P} (H)$, $A\in\Pi$.

The infinite-dimensional "Morse theory" of the action functional
gives rise to the {\it Floer homology} -- an analogue of the
classical Morse homology (for details
see e.g. \cite{Ho-Sa}, \cite{Sal}).

Namely, let $H$ be a
(time-dependent) Hamiltonian and let $J$ be an almost complex
structure on $M$ compatible with $\omega$.
For a generic pair $(H,J)$ -- which will be called
{\it a regular Floer pair} -- one can define a chain
complex $CF_\ast (H,J)$. It is a graded
module over the {\it anti-graded} ring $\Lambda_\omega$ (see
Section~\ref{subsec-Novikov-quantum}) equipped with a
$\Lambda_\omega{\hbox{\rm -linear}}$
differential that decreases the grading by 1. As a
$\Lambda_\omega{\hbox{\rm -module}}$ it is
generated by elements of ${\cal P} (H)$ and the module
structure fits with the action of $\Pi$ on ${\cal P} (H)$ defined by
the formula (\ref{eqn-cal-P-action-of-Pi-one-Ham}). Namely, if
$\hat{\gamma}\in {\cal P} (H)$ and $e^{ A}\in \Lambda_\omega$, then
the product $e^{A}\cdot \hat{\gamma}$
in the module $CF_\ast (H,J)$ over $\Lambda_\omega$ is equal to
$A\sharp\hat{\gamma}$.

Elements $\hat{\gamma}\in {\cal P} (H)$ are graded by the
Conley-Zehnder index $\mu({\hat{\gamma}})$ \cite{Co-Ze}. We follow
the sign convention which in the case of a small time-independent
Morse Hamiltonian $H$ gives the following relation: if $\gamma, f$
are constant maps to a critical point $x$ of $H$ and
$\mu_{Morse} (x)$ is the Morse index of $x$ then
\begin{equation}
\label{eqn-cz-index}
\mu ([\gamma ,f]) = 2n - \mu_{Morse} (x).
\end{equation}
The Conley-Zehnder index satisfies the property
\[
\mu (A\sharp \hat{\gamma}) = \mu (\hat{\gamma}) - 2 c_1 (A),
\]
so that the grading on $CF_\ast (H, J)$ fits indeed with the
structure of $CF_\ast (H, J)$ as a graded module over the {\it
anti-graded} ring $\Lambda_\omega$.

The homology $HF_\ast (H,J)$ of the chain complex $CF_\ast
(H,J)$ is called {\it the Floer homology group } -- it is a
graded module over the {\it anti-graded} ring $\Lambda_\omega$.
Using $\Lambda_\omega^0{\hbox{\rm -valued}}$ cochains one can
build an appropriate dual cochain complex $CF^\ast (H, J)$ whose
cohomology $HF^\ast (H, J)$ is called {\it the Floer cohomology
group} of $(M,\omega)$ -- it is a graded module over the {\it
graded} ring $\Lambda_\omega$. Similarly to the situation with
quantum (co)homology (see Section~\ref{subsec-Novikov-quantum})
there exists a non-degenerate evaluation pairing $(\cdot ,\cdot ):
HF^\ast (H, J) \times HF_\ast (H, J) \to \Lambda_\omega^0$ (see e.g.
\cite{PSS}), which, according to the universal coefficient
theorem, leads to a {\it group} isomorphism
\[
HF^k (H, J) = {\hbox{\rm Hom}} (HF_k (H, J), \Lambda_\omega^0 ),
\]
where ${\hbox{\rm Hom}} (HF_k (H, J), \Lambda_\omega^0 )$ denotes
the group of all {\it group homomorphisms} $HF_k (H, J)\to
\Lambda_\omega^0$. As in the classical Morse theory chains in
$CF_\ast (H,J)$ can be viewed as cochains in $CF^\ast (\bar{H},
J)$ leading to the Poincar{\'e} isomorphism of
$\Lambda_\omega^0{\hbox{\rm -modules}}$:
\[
PD: HF^\ast (H,J) \to HF_{2n-\ast} (\bar{H},J),
\]
which, as in the situation with quantum (co)homology (see
Section~\ref{subsec-Novikov-quantum}), can be extended to an
isomorphism of $\Lambda_\omega{\hbox{\rm -modules}}$.
As an additive module over $\Lambda_\omega$ the cohomology
$HF^\ast (H, J)$ is isomorphic to $H^\ast (M, \Q) \otimes_{\Q}
\Lambda_\omega$.

The Floer cohomology $HF^\ast (H, J)$ carries
a remarkable multiplicative structure called  {\it
the pair-of-pants product} (see \cite{PSS}, \cite{Sch-PhD}). This
multiplicative structure is isomorphic to the multiplicative
structure on the quantum cohomology of $(M,\omega)$.
Namely, S.Piunikhin, D.Salamon and
M.Schwarz
\cite{PSS} defined for
each regular Floer pair $(H,J)$ a
{\it $\Lambda_\omega{\hbox{\it -linear}}$
ring isomorphism preserving the grading}:
\[
\Psi_{H, J} : QH^\ast (M,\omega)\to HF^\ast (H, J).
\]
Similarly (see \cite{PSS}) there exists a
$\Lambda_\omega{\hbox{\rm -linear}}$ grading-preserving
isomorphism between quantum and Floer homology which together with
$\Psi_{H, J}$ intertwines the Poincar{\'e} isomorphism and evaluation
pairing in quantum and Floer (co)homology. The Floer (co)homology
groups for different regular Floer pairs $(H, J)$ are related by natural
isomorphisms which are compatible with the
Piunikhin-Salamon-Schwarz identifications of Floer and quantum
(co)ho\-mo\-lo\-gy.

%\bigskip
\subsection{Spectral numbers}
\label{subsect-spec-numbers-intro}

At first let $(H, J)$ be a regular Floer pair. Then
$(\bar{H}, J)$ is a regular Floer pair as well.
 Under our sign
convention -- see (\ref{eqn-def-ham-vf}) and
(\ref{eqn-def-action}) -- the Floer-Morse gradient-type connecting
trajectories, defining the differential in $CF_\ast (H, J)$,
correspond to the {\it downward} gradient flow of the action
functional ${\cal A}_H$ (and if $H$ is a sufficiently small
time-independent Morse function they also correspond to the {\it
upward} gradient flow of $H$).

Consider now a real filtration on $CF_\ast (H, J)$ defined by
${\cal A}_H$. Namely, let $CF_\ast^{(-\infty, s]} (H, J)$,
$-\infty\leq s\leq +\infty$, be spanned {\it over} $\Q$ by
all  $\hat{\gamma}\in {\cal P} (H)$ with
${\cal A}_H (\hat{\gamma}) \leq s$. Clearly, $CF_\ast^{( -\infty,
s]} (H, J)$ is invariant with respect to the differential of
$CF_\ast (H, J)$ and thus is a chain subcomplex. Note that $CF_\ast^{(
-\infty, s]} (H, J)$ is a (possibly infinite-dimensional) $\Q{\hbox{\rm
-subspace}}$ of $CF_\ast (H, J)$ and {\it not} a
$\Lambda_\omega{\hbox{\rm -submodule}}$ of $CF_\ast (H, J)$ since
multiplication by an element of $\Lambda_\omega$ may increase the
action.
Denote by $HF^s_\ast (H, J)$ the image of the homology of $CF_\ast^{(
-\infty, s]} (H, J)$ in $HF_\ast (H, J)$ under the inclusion
$CF_\ast^{( -\infty, s]} (H, J)\hookrightarrow CF_\ast (H, J)$.
Again note that $HF^s_\ast (H, J)$ is only
a $\Q{\hbox{\rm -subspace}}$ of
$HF_\ast (H, J)$ and {\it not} its $\Lambda_\omega{\hbox{\rm
-submodule}}$.

%\smallskip
\begin{defin}
\label{defin-spectral-number} {\rm Given a non-zero
class $\alpha\in QH^\ast (M, \omega)$ set
\[
c (\alpha , H) : = \inf\, \{\, s\, |\, PD (\Psi_{\bar{H}, J}
(\alpha)) \in HF^s_{2n-\ast} (H, J)\, \}.
\]
}
\end{defin}
%\bigskip

The idea behind the definition is due to C.Viterbo \cite{Vit} who
originally used it in the context of Morse homology and
finite-dimensional genera\-ting functions for Hamiltonian
symplectomorphisms of $\R^{2n}$. It has taken a considerable
work involving new techniques to extend this idea to the
Floer theory with various versions of Floer homology viewed as
infinite-dimensional analogs of Morse homology and with the action
functional viewed as an "infinite-dimensional generating
function". In the case of a cotangent bundle and the Floer theory
for Lagrangian intersections it was done by Y.-G. Oh \cite{Oh1},
\cite{Oh2}, then in the case of a closed symplectically aspherical
symplectic manifold and the Floer theory for Hamiltonian
symplectomorphisms by M.Schwarz \cite{Sch-1}, and in the case of a
general closed symplectic manifold and the Floer theory for
Hamiltonian symplectomorphisms it was done by Y.-G. Oh in his
recent papers \cite{Oh-act},\cite{Oh-new}, where
Definition~\ref{defin-spectral-number} is taken from.

The properties of spectral numbers that we need are
summarized in the following proposition.

%\smallskip
\begin{prop}
\label{prop-properties-c-general-case} Let
$\alpha\in QH^\ast (M, \omega)$ be non-zero.

\medskip
\noindent 1) If $H$ belongs to a regular Floer pair then $c
(\alpha, H)$ is finite. The number $c (\alpha, H)$
does not depend on the choice of $J$ in a regular Floer pair
$(H,J)$  and depends continuously on $H$ with respect to the
$C^0{\hbox{\it -norm}}$ on the space of (time-dependent)
Hamiltonian functions. Thus one can define $c (\alpha, H)$ for
any $H$ (i.e. not
necessarily belonging to a regular Floer pair). The resulting quantity
depends $C^0{\hbox{\it -continuously}}$ on $H$.

\medskip
\noindent 2) $c (\lambda e^A \alpha , H) = c (\alpha , H) - \omega
(A)$ for any $H$, $A\in\Pi$ and any non-zero
$\lambda\in \Q$.

\medskip
\noindent 3) If
$H^\prime = H + h (t)$ for some function $h:
S^1\to \R$ then $c (\alpha, H^\prime) = c (\alpha, H) -
\int_0^1 h (t) dt$.

\medskip
\noindent 4) Pick any
$\epsilon
> 0$ and any Floer chain representing
the Floer homology class
$PD (\Psi_{\bar{H}, J} (\alpha ))
\in HF_{2n-\ast} (H, J)$. Decompose the chain into a sum
$\sum_i \lambda_i \hat{\gamma}_i$,
$\lambda_i\in \Q$, $\hat{\gamma}_i\in {\cal P} (H)$. Then this sum
has to contain a non-zero term
$\lambda_i \hat{\gamma}_i$ with  ${\cal
A}_H (\hat{\gamma_i}) \geq c (\alpha, H) - \epsilon$. Moreover, one can
prove that $c (\alpha, H)$ belongs to the
spectrum ${\hbox{\it Spec}}\, (H)$.

\medskip
\noindent 5)If $\alpha\in H^\ast (M,\Q)\subset QH^\ast (M,\omega)$
is a {\bf singular}
cohomology class then
\begin{equation}
\label{eqn-c-spec-ineqs} - \int_0^1 \sup_M H^t dt \leq c (\alpha,
H)\leq - \int_0^1 \inf_M H^t dt.
\end{equation}

\medskip
\noindent 6) If a quantum product $\alpha_1\ast \alpha_2$ is
non-zero and Hamiltonians $H_1, H_2$ are normalized  then
\[
c (\alpha_1\ast \alpha_2, H_1\sharp H_2) \leq c (\alpha_1, H_1) +
c (\alpha_2, H_2).
\]

\medskip
\noindent 7) Assume that $H_1, H_2$ are normalized Hamiltonians
such that $\widetilde{\varphi}_{H_1} = \widetilde{\varphi}_{H_2}$
in $\tHam  (M,\omega)$. Then $c (\alpha, H_1)
= c (\alpha, H_2)$. Thus for any $\widetilde{\varphi}\in\tHam (M,\omega)$
one can define
\[
c (\alpha, \widetilde{\varphi}) = c (\alpha, H),
\]
where $H$ is any normalized Hamiltonian such that
$\widetilde{\varphi} = \widetilde{\varphi}_H$.

\medskip
\noindent 8) The function
$\widetilde{\varphi}\mapsto c (\alpha, \widetilde{\varphi})$ is
invariant under conjugation in $\tHam  (M,\omega)$.

\medskip
\noindent 9) If $(M, \omega)$ is symplectically aspherical then $c
({\bf 1}, \widetilde{\varphi}) = - c (m, \widetilde{\varphi}^{-1})$.
\end{prop}

The proof of the proposition will be outlined in
Section~\ref{sect-properties-spec-numbers-proofs}. The last equality
in part 9 is a partial case of a more general
statement relating spectral numbers of inverse elements
 in
$\tHam  (M,\omega)$ -- see \cite{EP}.

%\bigskip
\subsection{The main theorems}
\label{subsec-thm-action-length}

\begin{thm}
\label{thm-spectrum-length}

Let $(M,\omega)$ be a closed connected strongly semi-positive
symplectic manifold. Let $H$ be a normalized Hamiltonian and let
$E\in QH^\ast (M,\omega)$ be the Euler class. Assume
that for a positive integer $g$ the class $E^g\in QH^\ast
(M,\omega)$ is non-zero and  $c (E^g, H) >0$.
Then ${\hbox{\it cl}}\,(\widetilde{\varphi}_H) > g$ in
the group $\tHam  (M,\omega)$.
\end{thm}
%\smallskip

The proof of Theorem~\ref{thm-spectrum-length} occupies
Sections~\ref{sec-upsilon-K-area}-\ref{sect-pf-thm-spectrum-length-part-II}.
Its general overview and final steps can be found in
Section~\ref{sect-pf-thm-spectrum-length-part-II}.

In the case when $(M, \omega)$ is symplectically aspherical and
the Euler characteristic
$\chi (M)$ is zero the class $E$ vanishes and
Theorem~\ref{thm-spectrum-length} cannot be applied. However the
problem can be fixed by means of the following result.

%\smallskip
\begin{thm}[L.Polterovich]
\label{thm-sympl-aspher-case}

Let $(M,\omega)$ be a closed connected
symplectically aspherical manifold. Let $H$
be a normalized Hamiltonian. Assume that $c (m,\widetilde{\varphi}) > 0$. Then
${\hbox{\it cl}}\,(\widetilde{\varphi}) > 1$ in the group
$\tHam  (M,\omega)$.
\end{thm}
%\smallskip

The proof, unlike in the case of
Theorem~\ref{thm-spectrum-length}, involves only a short
calculation using basic properties of spectral numbers -- see
Section~\ref{subsect-pf-thm-sympl-aspher-case}.

Normalized Hamiltonian functions giving rise to positive spectral
numbers are produced by the following construction.

%\smallskip
\begin{prop}[Y.Ostrover, \cite{Ostr}]
\label{prop-ostrover}
Let $(M,\omega)$ be a closed connected strongly semi-posi\-tive
symplectic manifold.
Let $B\subset M$ be a ball and let $F: S^1\times M\to \R$ be a
normalized Hamiltonian such that $\varphi_F$ displaces $B$. Then
for any $\alpha\neq 0$
\[
c (\alpha, F\sharp w H_B) = c (\alpha, F) + w
\]
and thus $c (\alpha, F\sharp w
H_B) > 0$ for any sufficiently large $w$.
\end{prop}
%\smallskip

For a proof of Proposition~\ref{prop-ostrover} see
Section~\ref{subsect-pf-ostrover}.
Proposition~\ref{prop-ostrover} combined with
Theorems~\ref{thm-spectrum-length} yields
the following immediate corollary.

%\smallskip
\begin{cor}
\label{cor-infinite-cl-spectrum}

Let $(M,\omega)$ a be closed, connected and strongly semi-posi\-tive
symplectic manifold. Let $E\in QH^\ast (M,\omega)$ be the Euler
class.
Assume that the class $E^g\in QH^\ast (M,\omega)$ is non-zero for
any non-negative $g$. Then ${\hbox{\it
cl}}\,(\tHam  (M,\omega)) = +\infty$.

\end{cor}
%\smallskip

\begin{rem}
\label{rem-euler-class-ep}
{\rm
It is shown in \cite{EP} that if $(M, \omega)$ is spherically monotone and $E$
is {\it invertible} in the ring $QH^\ast (M,\omega)$ then there
exists
a non-trivial homogeneous quasimorphism on $\tHam  (M,\omega)$ and hence
${\hbox{\it
cl}}\,(\tHam  (M,\omega)) = +\infty$.
}
\end{rem}

%\bigskip
\subsection{Applications of
Theorem~\ref{thm-spectrum-length} }
\label{subsec-thm-action-length-applications}

We say that the class
$E^g\in QH^\ast (M,\omega) = H^\ast (M, \Q)\otimes_{\Q} \Lambda_\omega$
is {\it split} if it can be
represented as
\begin{equation}
\label{eqn-m-g-split} E^g = \alpha\otimes_{\Q} e^{C},
\end{equation}
$\alpha\in H^\ast (M, \Q)$, $e^{C}\in
\Lambda_\omega$.

%\smallskip
\begin{cor}
\label{cor-comm-length-quantum-length-ostrover}
Let $(M,\omega)$ a be closed, connected and strongly semi-posi\-tive
symplectic manifold.
Let $B\subset M$ be a ball and let $F: S^1\times M\to \R$ be a
normalized Hamiltonian such that $\varphi_F$ displaces $B$.

%\smallskip
\noindent 1) Assume that the class $E^g\in QH^\ast (M, \omega)$ is
non-zero and split as $E^g = \alpha\otimes_{\Q} e^{C}$ for a {\bf
particular} $g>0$. Then for any
\[
w > \int_0^1 \sup_M F (t, \cdot  ) dt + \omega (C)
\]
one has ${\hbox{\it cl}}\,(\widetilde{\varphi}_{F\sharp w H_B}) >
g$ in the group $\tHam  (M,\omega)$.

%\smallskip
\noindent 2) Assume that for {\bf each} $g>0$ the class $E^g$ is
non-zero and split as $E^g = \alpha_g\otimes_{\Q} e^{C_g}$. Set $I_g :=
\omega (C_g)$ and assume that $I:=\lim_{g\to\infty} I_g/g >0$.
Then for any $w>0$
\[
\| \widetilde{\varphi}_{wH_B}\|_{cl} \geq w/I.
\]

\end{cor}
%\smallskip

Below we present two examples, where the multiplicative structure of
the quantum cohomology ring of $(M,\omega)$ is known. In these
cases the classes $E^g$ are non-zero and split  and the quantities $I_g =
\omega (C_g)$ are computable.

%\bigskip
\begin{exam}
\label{exam-cp} {\rm Let $M = \C P^n$ and let
$\omega$ be the standard Fubini-Study form. Let $A, a$ be
respectively, the generators of $H_2 (\C P^n)$ and
$H^2 (\C P^n)$ which are positive with respect to $\omega$.
Set $q = e^{ A}\in \Lambda_\omega$. Then the
(rational) quantum cohomology ring of $(\C P^n, \omega)$ is
isomorphic to the graded polynomial ring
\[
\frac{\Q [a, q, q^{-1}]}{ \{ a^{n+1} = q,\  q\cdot q^{-1} =
1\} },
\]
where the degree of $a$ is two and the degrees of $q$ and $q^{-1}$
are, respectively, $2n+2$ and $- 2n -2$ \cite{Ru-Ti},
\cite{Ru-Ti-1}, \cite{Wi}.

Since the class $A$ generates $\Pi = H_2 (\C P^n)$ and $c_1
(A) = n + 1$ we can easily check that the Euler class $E$ actually
belongs to the singular cohomology $H^{2n} (\C P^n , \Q)$ and can be written as $E = \chi (\C P^n) m$. Now in the
quantum cohomology ring of $\C P^n$ one obviously has $\chi
(\C P^n ) m = \chi (\C P^n ) a^n$ and thus
\[
\chi (\C P^n ) m^g = \chi (\C P^n ) a^{gn} = \chi (\C P^n ) a^l\otimes_{\Q} e^{ k A}\in H^\ast (\C P^n , \Q )
\otimes_{\Q} \Lambda_\omega ,
\]
where $gn = k (n+1) + l$ with $k, l\in {\bf Z}, 0\leq l\leq n$.
Therefore $E^g = (\chi (\C P^n ) m)^g$ is split
for any $g$ and
\[
I_g = k\omega (A),
\]
where $k = [ gn / n+1 ]$. This yields
\begin{equation}
\label{eqn-I-g-I-CP-n}
I:= \lim_{g\to\infty} \frac{I_g}{g} = \lim_{g\to\infty} \frac{[ gn /
n+1 ]\, \omega (A)}{g} = \frac{n}{n+1} \omega (A).
\end{equation}

}
\end{exam}
%\smallskip

%\bigskip
\begin{exam}
\label{exam-grassm} {\rm Consider a generalization of the
previous example: let $M = Gr\, (r, n)$, $1\leq r\leq n-1$,
be the Grassmannian of complex $r{\hbox{\rm -dimensional}}$
subspaces in $\C^n$. The manifold $Gr\, (r, n)$ carries
a natural symplectic form $\omega$
--
one can view the Grassmannian as the result of symplectic
reduction for the Hamiltonian action of $U(r)$ (by multiplication
from the right) on the space of complex $r\times n$ matrices which
carries a natural symplectic structure.

The linear basis of the cohomology group $H^\ast (Gr\, (r, n),
\Q)$ is formed by the Schubert classes $\sigma^I$, where
$I=(i_1, \ldots, i_k)$ is a partition such that $r\geq i_1\geq
\ldots\geq i_k > 0$, $k\leq n-r$. Given two numbers $i,j$ $0\leq
i\leq r$, $0\leq j\leq n-r$, denote by $\{ i\times j\}$ the
partition $(i,i,\ldots , i)$ with $i$ repeated $j$ times. (If $i
=0$ or $j =0$ then $\{ i\times j\}$ is an empty partition). In
particular, the top-dimensional cohomology class $m\in H^{2 r
(n-r)} (Gr\, (r, n), \Q )$ can be viewed as $\sigma^{\{
r\times (n-r)\} }$.

Let $A$ be the generator of $H_2 (Gr\, (r, n), {\bf Z})
\cong {\bf Z}$ such that $\omega (A)>0$.
The quantum cohomology $QH^\ast (Gr\, (r, n),
\omega)$ is isomorphic, as a group, to $\Q [\sigma^I ]\otimes_{\Q}
\Q [q, q^{-1}]$ where $q = e^{ A}$ and $I$ runs over
the set of all partitions as above. The multiplicative structure
of the quantum cohomology ring of $(Gr\, (r,n),\omega )$ is
described in \cite{Bertr}, \cite{Sieb-Ti}, \cite{Wit-grasm}. In
particular, one can show that the Euler class for in the quantum
cohomology of $(Gr\, (r,n),\omega )$ belongs to the
singular cohomology \cite{Abrams},\cite{Bertram-Euler}:
\begin{equation}
\label{eqn-euler-class-grassm} E = \chi (Gr\, (r,n)) m.
\end{equation}
Moreover
there exists an algorithm
\cite{Bertr-CiF-Fult} which allows to compute explicitly
the structural coefficients of the quantum cohomology ring
$QH^\ast (Gr\, (r, n), \omega)$ with respect to the basis $\{
\sigma^I\}$. Using this algorithm A.Postnikov obtained the
following result \cite{Post}.

%\smallskip
\begin{prop}
\label{prop-post} Fix $g\geq 1$ and consider the quantum
cohomology class
$m^g\in QH^\ast (Gr\, (r,n))$.
At least one of the following two
statements about the numbers $r,n,g$ is always true:

%\medskip
\noindent (I) $gr = an +b$, $a = [gr/n]$, $0\leq b \leq r$.

%\smallskip
\noindent (II) $g(n-r)=cn+d$, $c=[g (n-r)/r]$, $0\leq d \leq n-r$.

%\smallskip
\noindent In the case (I)
\[
m^g = \sigma^{\{ b\times (n-r)\} }\otimes_{\Q} q^{a (n-r)}.
\]

%\smallskip
\noindent In the case (II)
\[
m^g = \sigma^{\{ r\times d\} }\otimes_{\Q} q^{cr}.
\]

\end{prop}
%\smallskip

Recalling that $c_1 (A) = n \neq 0$, $\chi (Gr\, (r,n)) > 0$
and applying (\ref{eqn-euler-class-grassm}) together with
Proposition~\ref{prop-post} we get that  $E^g = (\chi
(Gr\, (r,n)) m)^g$ is non-zero and split for any $g$ with
\[
I_g = a (n-r) \omega (A) = \biggl[ \frac{gr}{n}\biggr] (n-r)
\omega (A)
\]
in the case (I) and
\[
I_g = c r \omega (A) = \biggl[ \frac{g (n-r)}{n}\biggr] r \omega
(A)
\]
in the case (II). Thus
\begin{equation}
\label{eqn-I-g-I-Grassm}
I:=\lim_{g\to\infty} \frac{I_g}{g} = \frac{r (n-r)}{n} \cdot\omega
(A).
\end{equation}

}
\end{exam}
\bigskip

\subsection{Proofs of
Corollary~\ref{cor-comm-length-quantum-length-ostrover} and
Theorems~\ref{thm-cl-cp-n}-\ref{thm-cp-n-grassms-infinite-cl}}
\label{sect-pfs-thm-cl-cp-n-thm-stable-norm-cp-n}

\medskip
\noindent {\bf Proof of
Corollary~\ref{cor-comm-length-quantum-length-ostrover}.}

%\smallskip
\noindent Since $E^g = \alpha\otimes_{\Q} e^C$ parts 2 and 5 of
Proposition~\ref{prop-properties-c-general-case} yield
\[
c (E^g, F) = c (\alpha, F) - \omega (C)\geq
- \int_0^1 \sup_M F^t dt - \omega (C).
\]
But according to Proposition~\ref{prop-ostrover},
\[
c (\alpha, F\sharp wH_B) = c (\alpha, F) +w.
\]
Combining the inequality above
with the hypothesis
\[
w > \int_0^1 \sup_M F (t, \cdot  ) dt + \omega (C),
\]
we immediately get
\[
c (E^g, F\sharp wH_B) >0.
\]
Therefore, according to
Theorem~\ref{thm-spectrum-length},
${\hbox{\it cl}}\,(\widetilde{\varphi}_{F\sharp w H_B}) >
g$.
Part 1 of
Corollary~\ref{cor-comm-length-quantum-length-ostrover} is proven.

To show part 2 of
Corollary~\ref{cor-comm-length-quantum-length-ostrover} we
start with
the following simple observations. Recall that
\[
\widetilde{\varphi}_{F\sharp kwH_B} = \widetilde{\varphi}_{F}\cdot
\widetilde{\varphi}_{k w H_B}
\]
for any $k\in\Z$ and $w$. This, together with the definition of
commutator length, yields
\[
{\hbox{\it cl}}\, (\widetilde{\varphi}_{F}) = {\hbox{\it cl}}\,
(\widetilde{\varphi}_{F}^{-1}) \geq | {\hbox{\it cl}}\,
(\widetilde{\varphi}_{F\sharp k w H_B})
-
{\hbox{\it cl}}\, (\widetilde{\varphi}_{kwH_B}) |.
\]
Since $H_B$ is time-independent,
$\widetilde{\varphi}_{wH_B}^k = \widetilde{\varphi}_{k wH_B}$ for
any integer $k$. Combining all these observations, we see that
\[
{\hbox{\it cl}}\, (\widetilde{\varphi}_{F})\geq | {\hbox{\it
cl}}\, (\widetilde{\varphi}_{F\sharp k wH_B})
-
{\hbox{\it cl}}\, (\widetilde{\varphi}_{ wH_B}^k ) | \geq
{\hbox{\it cl}}\, (\widetilde{\varphi}_{F\sharp kw H_B})
-
{\hbox{\it cl}}\, (\widetilde{\varphi}_{ wH_B}^k )
\]
for any positive integer $k$ and therefore
\[
{\hbox{\it cl}}\, (\widetilde{\varphi}_{ wH_B}^k ) \geq {\hbox{\it
cl}}\, (\widetilde{\varphi}_{F\sharp k wH_B}) - {\hbox{\it cl}}\,
(\widetilde{\varphi}_{F}).
\]
Dividing the last inequality by $k$ and taking the limit as $k$
goes to $+\infty$ one gets
\begin{equation}
\label{eqn-ineq-cor-st-norm} \| \widetilde{\varphi}_{w H_B}\|_{cl}
= \lim_{k\to +\infty} \frac{ {\hbox{\it cl}}\,
(\widetilde{\varphi}_{F\sharp k w H_B})}{k}.
\end{equation}
Now pick an arbitrary small $\delta >0$. Then
\[
\lim_{g\to +\infty} (g\, I\,  (1+\delta) - I_g) = +\infty.
\]
Hence for any sufficiently large integer $g$
\[
g\, I\, (1+\delta) > \int_0^1 \sup_M F (t, \cdot  ) dt + I_g.
\]
Introducing a new integral parameter $k$ such that
\[
g = \biggl[ \frac{kw}{I\, (1+\delta)} \biggr]
\]
we see that for any sufficiently large integer $k$
\[
kw > \int_0^1 \sup_M F (t, \cdot  ) dt + I_{[kw/I (1+\delta)]}.
\]
Thus, according to the already proven part 1, for all sufficiently
large integers $k$
\[
{\hbox{\it cl}}\,(\widetilde{\varphi}_{F\sharp kw H_B}) > \biggl[
\frac{kw}{I (1+\delta)} \biggr].
\]
Together with (\ref{eqn-ineq-cor-st-norm}) this yields
\[
\| \widetilde{\varphi}_{w H_B} \|_{cl}\geq \lim_{k\to +\infty}
\frac{[kw/I (1+\delta)]}{k} = \frac{w}{I (1+\delta)}.
\]
This is true for any sufficiently small $\delta>0$ and thus
\[
\| \widetilde{\varphi}_{w H_B} \|_{cl}\geq w/I.
\]
Part 2 of the corollary is proven.
\b
%\smallskip

\bigskip
\noindent {\bf Proof of Theorem~\ref{thm-cl-cp-n}.}

%\smallskip
\noindent Theorem~\ref{thm-cl-cp-n} follows immediately from part
1 of Corollary~\ref{cor-comm-length-quantum-length-ostrover} and
formula (\ref{eqn-I-g-I-CP-n}) of
Example~\ref{exam-cp}.
\b
%\smallskip

\bigskip
\noindent {\bf Proof of Theorem~\ref{thm-stable-norm-cp-n}.}

%\smallskip
\noindent The theorem follows from part 2 of
Corollary~\ref{cor-comm-length-quantum-length-ostrover}
and Example~\ref{exam-cp}.
\b
%\smallskip

\bigskip
\noindent {\bf Proof of
Theorem~\ref{thm-ostrover-examples-sympl-aspher}.}

%\smallskip
\noindent
As in the proof of part 1 of
Corollary~\ref{cor-comm-length-quantum-length-ostrover}
we easily get that if
$w > \int_0^1 \sup_M F^t dt$ then $ c (m, F\sharp wH_B) >0$
and, according to Theorem~\ref{thm-sympl-aspher-case},
${\hbox{\it cl}}\,(\widetilde{\varphi}_{F\sharp kw H_B}) > 1$.
The theorem is proven.

Note that for a symplectically aspherical
$(M,\omega)$ the
quantum and the singular cohomology rings
coincide and $E = \chi (M) m$. Therefore if $\chi (M)\neq 0$
one can prove the theorem by applying part 1 of
Corollary~\ref{cor-comm-length-quantum-length-ostrover}
instead of Theorem~\ref{thm-sympl-aspher-case}.
\b
%\smallskip

\bigskip
\noindent {\bf Proof of
Theorem~\ref{thm-cp-n-grassms-infinite-cl}.}

%\smallskip
\noindent Theorem~\ref{thm-cp-n-grassms-infinite-cl} follows from
part 2 of Corollary~\ref{cor-comm-length-quantum-length-ostrover}
and formula (\ref{eqn-I-g-I-Grassm}) of Example~\ref{exam-grassm}.
\b
%\smallskip

%\bigskip
%\bigskip
\section{Properties of spectral numbers -- proofs}
\label{sect-properties-spec-numbers-proofs}

The goal of this section is to outline the proof of
Proposition~\ref{prop-properties-c-general-case}. At the end of
the section we will also prove Proposition~\ref{prop-ostrover}.

\subsection{Outline of the proof of
Proposition~\ref{prop-properties-c-general-case} (cf. \cite{EP},
\cite{Oh-new})} \label{subsect-pfs-prop-properties-c-general-case}

\medskip
\noindent {\bf Part 1.} Part 1 is proven in \cite{Oh-new}.
For symplectically aspherical manifolds it was first proven in
\cite{Sch-1}; for a singular cohomology class
$\alpha$ and a general symplectic manifold
it was proven
before in \cite{Oh-act}. A proof for spherically monotone
symplectic manifolds can be found in \cite{EP}.\ \b

\medskip
\noindent {\bf Part 2.} According to part 1, spectral numbers
depend continuously on the Hamiltonian. Therefore without loss of
generality we may assume that $H$ belongs to a regular Floer pair
$(H, J)$. Since $\Psi_{\bar{H}, J}$ is linear over $\Lambda_\omega$,
one can pass from Floer chains representing $PD (\Psi_{\bar{H}, J}
(\alpha))$ to Floer chains representing $PD (\Psi_{\bar{H}, J}
(\lambda e^{ A}\alpha))$ and backwards by multiplying them by
$\lambda e^{ A}$ or, respectively, by its inverse $\lambda^{-1} e^{-
A}$. According to  (\ref{eqn-action-funct-action-of-Pi}),
the multiplication by $\lambda e^{ A}$ decreases the level of a chain, with
respect to the
filtration on $CF_\ast (H, J)$, by $\omega (A)$ and the
multiplication by $\lambda^{-1} e^{-  A}$ increases it by the same
quantity yielding the necessary property.\ \b

\medskip
\noindent {\bf Part 3.} Follows immediately from
the definitions.\ \b

\medskip
\noindent {\bf Part 4.} Using the continuity of spectral
numbers with respect to Hamiltonian we may assume without loss of
generality that $H$ belongs to a regular Floer pair $(H, J)$.
Suppose the property does not hold. Then for some $\epsilon > 0$
the Floer homology class $PD (\Psi_{\bar{H}, J} (\alpha )) \in
HF_{2n-\ast} (H, J)$ can be represented by a chain of the form $\sum
b_i \hat{\gamma_i}$ (the sum may be infinite), with $b_i\in \Q$,
$\hat{\gamma}_i \in {\cal P} (H)$, so that ${\cal A}_H
(\hat{\gamma}_i) \leq c (\alpha, H) - \epsilon$ for all $i$. But
then $PD (\Psi_{\bar{H}, J} (\alpha )) \in HF^s_{2n-\ast} (H, J)$ for
$s= c (\alpha, H) -\epsilon$ which contradicts the definition of $c
(\alpha, H)$. Part 3 is proven. The statement
$c (\alpha, H)\in {\hbox{\it Spec}}\, (H)$ is proven in
\cite{Oh-new} (for a singular class $\alpha$ it was proven in
\cite{Oh-act}).\ \b

\medskip
\noindent {\bf Part 5.} This is proven in \cite{Sch-1} for
symplectically aspherical manifolds. The general case is proven in
\cite{Oh-act}.\ \b

\medskip
\noindent {\bf Part 6.} Part 6 is proven in \cite{Sch-1} (see
Proposition 4.1 there) for symplectically aspherical
manifolds (in that case the quantum cohomology coincides with the singular
one and the quantum product is the usual cup-product). In view
of the results from \cite{PSS} (also see \cite{ME}) the same proof works also for
an arbitrary strongly semi-positive manifold.
\ \b

\medskip
\noindent {\bf Part 7.} Part 7 is proven in the general case by
the same argument that was used in \cite{Sch-1} in the
symplectically aspherical case. Namely let $\{
H_\tau\}_{0\leq \tau\leq 1}$ be a homotopy $H_1$ and $H_2$ in the class
of normalized Hamiltonians generating $\widetilde{\varphi}_{H_1} =
\widetilde{\varphi}_{H_2}$. Similarly to
the situation in the symplectically aspherical case (see Proposition 3.1
in \cite{Sch-1}), there exists a natural
one-to-one correspondence between ${\cal P} (H_1)$ and ${\cal P}
(H_\tau)$ for any $\tau$ and hence a one-to-one correspondence
between ${\hbox{\it Spec}}\, (H_1)$ and ${\hbox{\it Spec}}\, (H_\tau)$.
Moreover, a computation similar to
Lemma 3.3 in \cite{Sch-1} shows that
${\hbox{\it Spec}}\, (H_1) = {\hbox{\it Spec}}\, (H_\tau)$ for any
$\tau$. According to Proposition~\ref{prop-spectr-measure-zero},
the set ${\hbox{\it Spec}}\, (H_1)\subset \R$ is of measure
zero. Now part 4 of the proposition says that for any $\tau$ the number $c
(\alpha, H_\tau)$ belongs to the same set ${\hbox{\it Spec}}\,
(H_1)$ of measure zero. On the other hand, according to part 1, $c (\alpha,
H_\tau)$ depends continuously on $\tau$. Therefore it is a
constant function of $\tau$ and  $c (\alpha, H_1) = c (\alpha,
H_2)$, which proves part 7 of the proposition.\ \b

\medskip
\noindent {\bf Part 8.} Take two normalized Hamiltonians $H_1$ and
$H_2$. In view of Part 7, it suffices to prove that $c (\alpha,
H_1) = c (\alpha, H_2\sharp H_1\sharp \bar{H_2})$. Set
\[
H^\prime := H_2\sharp H_1 \sharp \bar{H}_2,
\]
\[
H^{\prime\prime}: = H_1 (t, \varphi^{-1}_{H_2} (x)).
\]
Then $\widetilde{\varphi}_{H^\prime} =
\widetilde{\varphi}_{H^{\prime\prime}}$. Also, if $H_1, H_2$ are
normalized then $H^\prime$ and $H^{\prime\prime}$ are normalized
as well. Therefore, in view of part 7, one has $c (\alpha,
H^\prime ) = c (\alpha, H^{\prime\prime})$. On the other hand,
there exists a natural isomorphism between the
Floer complexes of $H_1$ and $H^{\prime\prime}$ preserving the filtration.
Hence $c (\alpha, H_1) = c (\alpha,
H^{\prime\prime})$ and therefore $c (\alpha, H_1) = c (\alpha,
H^\prime )$. This finishes the proof of part 8. \ \b

\medskip
\noindent {\bf Part 9.} This is proven in \cite{Sch-1}.\ \b

%\bigskip
\subsection{Proof of
Theorem~\ref{thm-sympl-aspher-case}}
\label{subsect-pf-thm-sympl-aspher-case}

Part 6 of
Proposition~\ref{prop-properties-c-general-case} yields
\begin{equation}
\label{eqn-thm-polt-1} c (m, a b a^{-1} b^{-1}) \leq c (m, a) + c ({\bf 1},
ba^{-1}b^{-1}).
\end{equation}
On the other hand, in view of parts 8 and 9 of
Proposition~\ref{prop-properties-c-general-case},
\begin{equation}
\label{eqn-thm-polt-2} c ({\bf 1}, ba^{-1}b^{-1}) = c ({\bf 1},
a^{-1}) = - c (m, a).
\end{equation}
Combining (\ref{eqn-thm-polt-1}) and (\ref{eqn-thm-polt-2}) we get
\begin{equation}
\label{eqn-thm-polt-3} c (m, a b a^{-1} b^{-1})\leq 0.
\end{equation}
Hence if $c (m, \widetilde{\varphi}) > 0$ then $\hbox{\it cl}\,
(\widetilde{\varphi}) > 1$ in
$\tHam  (M,\omega)$
and the theorem is proven.
\b

%\bigskip
\subsection{Proof of
Proposition~\ref{prop-ostrover}.} \label{subsect-pf-ostrover}

Since $\varphi_F$ displaces $B$ and $\varphi_{w H_B}$ is supported
in $B$, the fixed points of
$\widetilde{\varphi}_{F\sharp w H_B}$ coincide with the fixed points of
$F$. An easy calculation shows that if $\hat{\gamma}\in {\cal P} (F)$
then the value of ${\cal A}_{F\sharp wH_B}$ on
the element of
${\cal P} (F\sharp wH_B)$, corresponding to $\hat{\gamma}$, is equal to
${\cal A}_F (\hat{\gamma})+ w$.
Hence
\begin{equation}
\label{eqn-ostr-pf-1} {\hbox{\it Spec}}\, (F\sharp wH) =
{\hbox{\it Spec}}\, (F) + w.
\end{equation}

Parts 1 and 4 of Proposition~\ref{prop-properties-c-general-case}
say that
$c (\alpha, F\sharp wH_B)\subset {\hbox{\it
Spec}}\, (F\sharp wH)$  and  $c (\alpha, F\sharp wH_B)$
changes continuously with $w$.
But, according to Proposition~\ref{prop-spectr-measure-zero}, for any
$w$ the set ${\hbox{\it Spec}}\, (F\sharp wH)$ is of measure zero.
Together
with (\ref{eqn-ostr-pf-1}) this yields
\[
c (\alpha, F\sharp wH_B) = c (\alpha, F\sharp H_B) + w.
\]
The proposition is proven. \ \b
%\bigskip

%\bigskip
%\bigskip
\section{A function
$\Upsilon_{l,g}$ on sets of conjugacy classes and K-area}
\label{sec-upsilon-K-area}

In this section we start to prepare the tools needed for the
proof of
Theorem~\ref{thm-spectrum-length}.
Our setup throughout the section will be an arbitrary connected Lie group
$G$
equipped with a bi-invariant Finsler pseudo-metric. We will
define a function on a set of conjugacy classes. This function
will be
used to measure the distance from a given element of the group
(or, more precisely, from its conjugacy class) to the set of
elements whose commutator length is bounded by a certain fixed number.
Then we will show that the function can be described in
the language of connections on the trivial $G{\hbox{\rm -bundle}}$
over an oriented compact surface with boundary.

Our main example will be the group
$\tHam  (M,\omega)$. The pseudo-metric
on this group is defined as the lift of the Hofer metric from
$\Ham\,  (M,\omega)$.

%\bigskip
\subsection{A function
$\Upsilon_{l,g}$ on sets of conjugacy classes}
\label{subsec-upsilon-defs}

We are going to define a function $\Upsilon_{l,g}$ on sets of $l$
conjugacy classes by modifying the definition of a similar
function $\Upsilon_l$ in \cite{ME}.

Namely, given a connected Lie group $G$ and a positive integer $k$
denote by $G^{(k)}\subset  G$ the subset formed by all products of
$k$ commutators in $ G$. Set $ G^{(0)} = Id$. Each $ G^{(k)}$ is
invariant under conjugation and contains inverses of all its
elements. Obviously,
\[
\{ Id\} \subset  G^{(1)}\subset  G^{(2)}\subset\ldots \subset [ G,
G],
\]
\[
\bigcup_{k\geq 1}  G^{(k)} = [ G, G]
\]
and
\[
{\hbox{\it cl}}\,(x) = \min\, \{k\, |\, x\in  G^{(k) }\}
\]
for $x\in [ G, G]$.

Now suppose the Lie algebra
${\mathfrak g}$ of $G$ admits a bi-invariant norm. Such a norm defines
a bi-invariant Finsler norm $\|
\cdot \|$ on $T G$ (which smoothly depends on a point in
$G$). This allows us to measure lengths of paths in $G$: if
$\gamma: [0,1]\to G$ is a smooth path then
\[
{\hbox{\rm length}}\, (\gamma) := \int_0^1 \| \gamma^\prime (t)\|
dt.
\]
Thus one can define
a {\it bi-invariant pseudo-metric} $\rho$ on $G$ called {\it
Finsler pseudo-metric}: if $\phi, \psi\in G$ then
\[
\rho (\phi, \psi) := \inf_\gamma {\hbox{\rm length}}\, (\gamma)
\]
where the infimum is taken over all paths $\gamma$ connecting
$\phi$ and $\psi$. It is a pseudo-metric in the sense that it
satisfies the same conditions as a genuine metric (it is
non-negative, symmetric and satisfies the triangle inequality) but
may be
degenerate:
$\rho (x, y)$ may be zero if $x\neq y$.

%\smallskip
\begin{defin}
\label{def-upsilon-pos-genus} {\rm Given a bi-invariant Finsler
pseudo-metric $\rho$ on $G$ and a tuple ${\cal C} = ({\cal
C}_1,\ldots , {\cal C}_l)$ of conjugacy classes in $G$ define
\[
\Upsilon_{l,g}\, ({\cal C}) = \inf_{ x\in  G^{(g)}, \varphi_i\in
{\cal C}_i } \rho\, (Id, x\cdot\prod_{i=1}^l \varphi_i ).
\]
In particular, in the notation of \cite{ME} one has
$\Upsilon_{l,0}\, ({\cal C}) = \Upsilon_l\, ({\cal C})$, where
$\Upsilon_l\, ({\cal C})$ is defined as: $\Upsilon_l\, ({\cal
C}):= \inf_{\varphi_i\in {\cal C}_i } \rho\, (Id, \prod_{i=1}^l
\varphi_i )$. (Obviously, $\Upsilon_{l,g}\, ({\cal C})$ does not
depend on the order of classes ${\cal C}_i$ in ${\cal C}$).

}
\end{defin}
%\bigskip

Now consider a particular case when $l=1$ and observe that
conjugate group elements have the same commutator length.
Definition~\ref{def-upsilon-pos-genus} yields

%\smallskip
\begin{prop}
\label{prop-comm-length-upsilon} Let ${\cal C}$ be a conjugacy
class in $G$ that lies inside $[G,G]$. Then the following
conditions are equivalent:

%\medskip
\noindent (i) ${\hbox{\it cl}}\, (x) = g$ for any $x\in {\cal C}$;

%\smallskip
\noindent (ii) $\Upsilon_{1, k} ({\cal C}) = 0$ for $k\geq g$ and
$\Upsilon_{1, k} ({\cal C}) > 0$ for $1\leq k < g$.

\end{prop}
%\smallskip

%\bigskip
\subsection{K-area and
$\Upsilon_{l,g}$} \label{subsect-def-K-area}

We will now present an alternative way to describe the function
$\Upsilon_{l,g}$.

Let $\Sigma$ be a connected oriented Riemann surface of genus $g$
with $l\geq 1$ infinite cylindrical ends $\Sigma_1, \ldots
,\Sigma_l$. Fix an orientation-preserving identification
$\Phi_i:  [0,+\infty)\times S^1\to \Sigma_i$, $1\leq i\leq l$, of
each end $\Sigma_i\subset \Sigma$ with the standard oriented
cylinder $[0,+\infty) \times S^1$.

Fix a volume form $\Omega$ on $\Sigma$ so that $\int_\Sigma \Omega
= 1$. According to Moser's theorem \cite{Mos}, any two such volume
forms coinciding at infinity can be mapped into each other by a
compactly supported
diffeomorphism of $\Sigma$.

Let $G$ be a connected Lie group whose tangent bundle is equipped
with a bi-invariant Finsler norm that defines a bi-invariant
Finsler pseudo-metric $\rho$ on $G$ (see
Section~\ref{subsec-upsilon-defs}). Identify the Lie algebra
${\mathfrak g}$ of $G$ with the space of right-invariant vector
fields on $G$. Assume that $G$ acts effectively
on a connected manifold $F$, i.e. $G\to
{\hbox{\it Diff}}\, (F)$ is a monomorphism.

Below we will not use any results from the Lie
theory and thus all our considerations will hold even if $G$ is an
infinite-dimensional Lie group, like $G=
\tHam\,  (M,\omega)$ which will be our main example.

Consider the trivial $G{\hbox{\rm -bundle}}$ $\pi: P\to\Sigma$,
with the fiber $F$.  Let $L^\nabla$
denote the curvature of a connection $\nabla$ on the bundle $\pi:
P\to\Sigma$.\footnote{We work only with $G{\hbox{\rm -connections}}$,
i.e. the connections whose parallel
transports belong to the structural group $G$.} If the fiber $\pi^{-1} (x)$ is identified with $F$
then to a pair of vectors $v,w\in T_x \Sigma$ the curvature tensor
associates an element $L^\nabla (v,w)\in {\mathfrak g}$. Here we
use the fact that $G$ acts effectively on $F$. If no
identification of $\pi^{-1} (x)$ with $F$ is fixed then $L^\nabla
(v,w)\in {\mathfrak g}$ is defined up to the adjoint action of $G$
on ${\mathfrak g}$. Thus if the norm $\|\cdot \|$ on $T G$
is the bi-invariant,
$\| L^\nabla (v,w) \|$ does not depend on the
identification of $\pi^{-1} (x)$ with $F$.

%\smallskip
\begin{defin}
\label{def-norm-curvature} {\rm Define $\|L^\nabla\|$ as
\[
{\|L^\nabla\|} = \max_{v,w} \frac{\| L^\nabla (v,w)\| }{|\Omega
(v,w) |},
\]
where the maximum is taken over all pairs $(v,w)\in T \Sigma
\times T \Sigma$ such that $\Omega (v,w)\neq 0$. }
\end{defin}
%\bigskip

Given a connection $\nabla$ on $\pi: P\to\Sigma$
its holonomy  along a loop based at $x\in\Sigma$ can be
viewed as an element of $G$ acting on $F$ provided that the bundle
is trivialized over $x$. If the trivialization of the bundle is
allowed to vary then the holonomy is defined up to conjugation in
the group $G$. Observe that the action of the gauge group does not
change $\|L^\nabla\|$.

Now let ${\cal C}=({\cal C}_1,\ldots,{\cal C}_l)$ be
conjugacy classes in $G$.

%\smallskip
\begin{defin}
\label{def-cal-L} {\rm Let ${\cal L} ({\cal C})$ denote the set of
connections $\nabla$ on $\pi: P\to \Sigma$ which are flat over the
set $\bigcup_{i=1}^l \Phi_i ([K,+\infty)\times S^1 )$ for some
$K>0$ and such that for any $i=1,\ldots,l$ and any $s\geq K$ the
holonomy of $\nabla$ along the oriented path $\Phi_i (s\times
S^1)\subset \Sigma_i$, with respect to a fixed trivialization of $\pi$,
does not depend on $s$ and lies in ${\cal
C}_i$. }
\end{defin}
%\bigskip

%\smallskip
\begin{defin}
\label{def-K-area} {\rm The number $0< {\hbox{\it K-area}}_{l,g}\,
({\cal C})\leq +\infty$ is defined as
\[
{\hbox{\it K-area}}_{l,g}\, ({\cal C}) = \sup_{\nabla\in {\cal L}
({\cal C}) }\|L^\nabla\|^{-1}.
\]

}
\end{defin}
%\bigskip

Obviously, since none of the ends of $\Sigma$ is preferred over
the others, the quantity ${\hbox{\it K-area}}_{l,g}\, ({\cal
C}_1,\ldots, {\cal C}_l)$ does not depend on the order of the
conjugacy classes ${\cal C}_1,\ldots, {\cal C}_l$.

Now let $\Upsilon_{l,g}$ be defined by means of the same
pseudo-metric on $G$ which was used in the definition of K-area.

%\smallskip
\begin{prop}
\label{prop-upsilon-k-area-pos-genus}
\[
\Upsilon_{l,g}\, ({\cal C})
=
\frac{1}{{\hbox{\it K-area}}_{\, l,g}\, ({\cal C})}.
\]
(If the K-area is infinite its inverse is assumed to be zero).
\end{prop}
%\smallskip

%\bigskip
%\bigskip
\subsection{Proof of
Proposition~\ref{prop-upsilon-k-area-pos-genus}}
\label{subsect-pf-k-area-upsilon}

As in the case of genus zero (see \cite{ME}), one needs to deal
with {\it certain systems of paths}.

%\smallskip
\begin{defin}
\label{def-system-of-paths-pos-genus} \rm{ An {\it
$(l,g){\hbox{\it -system}}$ of paths } is a tuple
\[
(a,b) = (a_1,\ldots, a_l, b_1, {\bar b}_1,\ldots , b_g, {\bar b}_g
)
\]
of $l+2g$ smooth paths $a_1,\ldots, a_l, b_1, {\bar
b}_1,\ldots , b_g, {\bar b}_g: [0,1]\to G$ such that

\medskip
\noindent 1) $\prod_{i=1}^l a_i (0)\cdot \prod_{k=1}^g b_k
(0)\cdot {\bar b}_k (0) = Id$;

\medskip
\noindent 2) $b_k (1)$ and ${\bar b}_k^{-1} (1)$ are conjugate for
each $k=1,\ldots, g$. }
\end{defin}
%\bigskip

The ${\hbox{\it length}}\, (a,b)$ of an $(l,g){\hbox{\rm
-system}}$ of paths $(a,b)$ is defined as the sum of lengths of
the paths that form the system, where the length of a path is
measured with respect to the chosen Finsler pseudo-metric on $G$.

%\smallskip
\begin{defin}
\label{def-system-of-paths-conj-classes-pos-genus} \rm{ Let ${\cal
C} = ({\cal C}_1,\ldots,{\cal C}_l)$ be conjugacy classes in
$G$. Define ${\cal G}_{l,g}\, ({\cal C})$ as the set of all
$(l,g){\hbox{\rm -systems}}$ of paths $(a,b)$ such that $a_i
(1)\in {\cal C}_i$, $i=1,\ldots,l$. }
\end{defin}
%\bigskip

The proof of Proposition~\ref{prop-upsilon-k-area-pos-genus}
follows immediately from the next two
propositions:

%\smallskip
\begin{prop}
\label{prop-k-area-via-systems-of-paths-pos-genus}
\[
\frac{1}{{\hbox{\it K-area}}_{\, l,g}\, ({\cal C})}
=
\inf_{\scriptscriptstyle (a,b)\in {\cal G}_{l,g}\, ({\cal C})}
{\hbox {\it length}}\, (a,b).
\]
\end{prop}
%\smallskip

%\smallskip
\begin{prop}
\label{prop-upsilon-via-systems-of-paths-pos-genus}
\[
\Upsilon_{l,g}\, ({\cal C}) = \inf_{\scriptscriptstyle (a,b)\in
{\cal G}_{l,g}\, ({\cal C})} {\hbox {\it length}}\, (a,b).
\]
\end{prop}
%\smallskip

Proposition~\ref{prop-upsilon-via-systems-of-paths-pos-genus}
extends Proposition 3.5.3 in \cite{ME}. We will prove it below.

Proposition~\ref{prop-k-area-via-systems-of-paths-pos-genus} can
be proved by exactly the same methods as in Section 10 of
\cite{ME}. Here we will discuss only the most important point
which needs to be understood as one carries over the proof from
\cite{ME} to the current situation.

Without loss of generality one can view $\Sigma$ as a compact
surface with boundary (rescale the
infinite cylindrical ends of $\Sigma$ to make them finite and then
take the closure of the open surface). The $(l,g){\hbox{\rm
-system}}$ of paths associated to a connection appears if we cut
open each of the $g$ handles of $\Sigma$ to get a surface
$\Sigma^\prime$ of genus zero with $l+2g$ boundary components.
Then, in the same way as in Section 10 of \cite{ME}, we construct
$l+2g$ \break paths $(a_1,\ldots, a_l, b_1, {\bar b}_1,\ldots ,
b_g, {\bar b}_g )$, with each path corresponding to a
boun\-dary component of $\Sigma^\prime$, so that
\[
a_1 (0)\cdot\ldots\cdot a_l (0)\cdot b_1 (0)\cdot \bar{b}_1
(0)\cdot\ldots\cdot b_g (0)\cdot \bar{b}_g (0) = Id\]
(cf.
Definition 3.5.1 in \cite{ME}). In particular, each
pair  $b_k, {\bar b}_k$, $k=1,\ldots ,g$, corresponds to
the two boundary components $T_k, {\overline T}_k$ of
$\Sigma^\prime$ that come from the cut along a circle in the
$k{\hbox{\rm -th}}$ handle of $\Sigma$.
The fact that each $b_k (1)$ is conjugate to ${\bar
b}_k^{-1} (1)$ comes from the difference of the trivializations of
the principal bundle $\Sigma\times G\to \Sigma$ over $T_k$ and
${\overline T}_k$. Indeed, recall that in the construction of the
paths $(a_1,\ldots, a_l, b_1, {\bar b}_1,\ldots ,
b_g, {\bar b}_g )$
(see Section 11.2.2 of \cite{ME}) the trivializations
of $\Sigma\times G\to \Sigma$ near the boundary components are not
arbitrary. However -- see Section 11.2.2 of \cite{ME}
--
as one identifies the oriented circles $T_k\cong {\overline
T}_k^{-1}$, $k=1,\ldots ,g$, one can always assume that the gauge
transformation that identifies the trivializations over $T_k$ and
${\overline T}_k^{-1}$ is a {\it constant} map $T_k\to G$. For
an appropriate global trivialization of the bundle $\Sigma\times
G\to \Sigma$ the constant image of $T_k\to G$ is simply a holonomy
of the connection over the closed path $b_k\circ {\bar b}_k^{-1}$
(the composition sign stands for the composition in the space of
paths) that goes once around the $k{\hbox{\rm -th}}$ handle of
$\Sigma$. Thus $(a_1,\ldots, a_l, b_1, {\bar b}_1,\ldots , b_g,
{\bar b}_g )$ is indeed an $(l,g){\hbox{\rm -system}}$ of paths.
The rest of the argument from \cite{ME} can be carried over to our
case in a straightforward manner. This finishes our discussion of
the proof of
Proposition~\ref{prop-k-area-via-systems-of-paths-pos-genus}.
\b

%\bigskip
\bigskip
\noindent {\bf Proof of
Proposition~\ref{prop-upsilon-via-systems-of-paths-pos-genus}.}

%\medskip
\noindent Pick any $g$ commutators $X_1,\ldots, X_g$ in $G$, where
$X_k = Y_k Z_k Y_k^{-1} Z_k^{-1}$, $k=1,\ldots, g$. Consider all
$(l,g){\hbox{\rm -systems}}$ of paths $(a,b) \in {\cal
G}_{l,g}\, ({\cal C})$ in which all paths, except for $a_1$,
are constant so that:

\begin{itemize}
\item{the image of each constant path
$a_i$, $i=2,\ldots, l$, is an element
$\phi_i\in {\cal C}_i$;}

\item{for each $k=1,\ldots , g$
the image of the constant path
$b_k$ is
$Y_k$ and the image of the constant path
${\bar b}_k$ is $Z_k Y_k^{-1} Z_k^{-1}$.}
\end{itemize}

Set $a_1 (1) = \phi_1\in {\cal C}_1$. Any such $(l,g){\hbox{\rm
-system}}$ of paths $(a,b)$ has to satisfy the condition
\[
Id = \bigg( \prod_{i=1}^l a_i (0) \bigg) \bigg( \prod_{k=1}^g b_k
(0) {\bar b}_k (0) \bigg)
 =
a_1 (0)\cdot \prod_{i=2}^l \phi_i \cdot \prod_{k=1}^g
Y_k Z_k Y_k^{-1} Z_k^{-1}
\]
and thus
\[
{\hbox{\it length}}\, (a,b)  = {\hbox{\it length}}\, (a_1) =   \]
\[ =
\rho\, (Id, a_1 (1)\cdot a_1^{-1} (0)) = \rho\, (Id, \bigg(
\prod_{i=1}^l \phi_i \bigg) \bigg( \prod_{k=1}^g X_k \bigg) ).
\]
Therefore
\[
\inf_{ (a,b)\in {\cal G}_{l,g}\, ({\cal C})} {\hbox{\it length}}\,
(a,b)\leq \inf_{\phi\in {\cal C}, \psi\in G^{(g)}} \rho\, (Id,
\bigg( \prod_{i=1}^l \phi_i \bigg)\cdot \psi )=
\Upsilon_{l,g}\, ({\cal C}).
\]

%\bigskip
Let us now prove the opposite inequality. Pick an arbitrary
$(a,b)\in {\cal G}_{l,g}\, ({\cal C})$ so that ${\bar b}_k (1) =
Z_k\cdot b_k^{-1} (1)\cdot Z_k^{-1}$ for some $Z_k\in G$,
$k=1,\ldots, g$. Then
\begin{eqnarray}
\label{eqn-upsilon-l-g-leq-syst-paths} {\hbox{\it length}}\, (a,b)
& =  & \sum_{i=1}^l \rho\, (Id, a_i (0)\cdot a_i^{-1} (1)) +
\nonumber\\ & + & \sum_{k=1}^g \rho\, (Id, b_k (0)\cdot b_k^{-1}
(1) ) \ +
\\
& + & \sum_{k=1}^g \rho\, (Id, Z_k\cdot b_k (1)\cdot Z_k^{-1}\cdot
{\bar b}_k (0) ). \nonumber
\end{eqnarray}
The right-hand side
of (\ref{eqn-upsilon-l-g-leq-syst-paths}) can be estimated from below
by means of the triangular inequality as
\begin{equation}
\label{eqn-upsilon-l-g-leq-syst-paths-1} {\hbox{\it length}}\,
(a,b)\geq \rho\, (Id,\, \Xi),
\end{equation}
where
\[
\Xi = \bigg( \prod_{i=1}^l a_i (0)\cdot a_i^{-1} (1) \bigg) \bigg(
\prod_{k=1}^g b_k (0)\cdot b_k^{-1} (1)\cdot Z_k\cdot b_k (1)\cdot
Z_k^{-1}\cdot {\bar b}_k (0) \bigg).
\]
Set $A_i = \prod_{j=1}^i a_j (0)$ and observe that
\[
\prod_{i=1}^l a_i (0)\cdot a_i^{-1} (1) = \bigg( \prod_{i=1}^l
A_i a_i^{-1} (1) A_i^{-1} \bigg) \cdot A_l = \bigg( \prod_{i=1}^l
\phi_i^{-1} \bigg) \cdot A_l,
\]
for some $\phi = (\phi_1,\ldots, \phi_l)\in{\cal C}$.
Denote by $X_k$ the commutator
$X_k = b_k^{-1} (1)\cdot Z_k\cdot b_k (1)\cdot Z_k^{-1}$. One can write
\begin{equation}
\label{eqn-upsilon-l-g-leq-syst-paths-2} \Xi = \bigg(
\prod_{i=1}^l \phi_i^{-1} \bigg) \cdot\Xi_1,
\end{equation}
where
\begin{equation}
\label{eqn-Xi-1}
\Xi_1 = A_l\cdot \prod_{k=1}^g b_k (0) X_k {\bar b}_k (0).
\end{equation}
Note that by the definition of an $(l,g){\hbox{\rm -system}}$
of paths
\[
A_l = {\bar b}_g^{-1} (0)\cdot b_g^{-1} (0)\cdot\ldots\cdot
{\bar b}_1^{-1} (0)\cdot b_1^{-1} (0).
\]
Also note that for any commutator $X$ and any $b\in G$ one has $X
b = b X^\prime$ for some other commutator $X^\prime$. Using these
observations we can move the commutators in the product in
(\ref{eqn-Xi-1}),
possibly changing them into other commutators, and rewrite $\Xi_1$ as:
\[
\Xi_1 = A_l\cdot \bigg( \prod_{k=1}^g X_k^\prime \bigg) \cdot
A_l^{-1}
\]
for some commutators $X_k^\prime$, $k=1,\ldots, g$. Thus
$\Xi_1\in G^{(g)}$. Now substituting
(\ref{eqn-upsilon-l-g-leq-syst-paths-2}) into
(\ref{eqn-upsilon-l-g-leq-syst-paths-1}) one gets
\[
{\hbox{\it length}}\, (a,b)\geq \rho\, (Id,\, \bigg( \prod_{i=1}^l
\phi_i^{-1} \bigg) \cdot \Xi_1),
\]
where $\phi\in {\cal C}, \Xi_1\in G^{(g)}$ and therefore
\[
\Upsilon_{l,g}\, ({\cal C}) = \inf_{\phi\in {\cal C}, \psi\in
G^{(g)}} \rho\, (Id, \bigg( \prod_{i=1}^l \phi_i^{-1} \bigg) \cdot
\psi ) \leq \inf_{(a,b)\in {\cal G}_{l,g}\, ({\cal C})} {\it
length}\, (a,b),
\]
which finishes the proof of the proposition.\ \b

%\bigskip
%\bigskip
\section{Hofer metric on the group
$\Ham\,  (M,\omega)$, Hamil\-tonian con\-nec\-tions and weak
coup\-ling } \label{sec-Ham-Hofer-metric-def}

We are going to apply the constructions from the previous section to
the group $G=\tHam  (M,\omega)$, where $(M,\omega)$ is
a closed connected symplectic manifold. In this case there exists
an additional construction -- {\it weak coupling} -- that allows
to estimate K-area from above and $\Upsilon_{l,g}$ from below.

{\sl Throughout this section we assume that $(M,\omega)$ is an
arbitrary closed connected symplectic manifold.}

The group $G=\tHam  (M,\omega)$ can be viewed as an
infinite-dimensional Lie group: it can be equipped with the
structure of an infinite-dimensional manifold so that the group
product and taking the inverse of an element become smooth
operations \cite{Ra-Schm}. The Lie algebra $\g$ of $\tHam
(M,\omega)$, which is also the Lie algebra of $\Ham\,  (M,\omega)$,
can be
identified with the Poisson-Lie algebra of all functions on $M$
with the zero mean value.

The norm $\| H \| = \max_M | H | $ on $\g$ is
bi-invariant and defines a bi-invariant
Finsler pseudo-metric on $G=\tHam  (M,\omega)$ which we will call
{\it Hofer pseudo-metric}. It is a lift of a famous
genuine metric on $\Ham\,  (M,\omega)$, called
{\it Hofer metric}
\cite{Ho},
\cite{L-McD-Hm},\cite{Pol-Hm}.
(The original metric introduced by Hofer was
defined by the norm $\| H\| = \max_M H - \min_M H$ and is
equivalent to the metric we use).

Let $\Sigma$ be a connected oriented Riemann
surface of genus $g$ with $l\geq 1$ infinite cylindrical ends and
a fixed area form $\Omega$ of total area 1 as
in Section~\ref{subsect-def-K-area}.
Consider the trivial bundle $\Sigma\times M\to
\Sigma$ with the fiber $F = (M,\omega )$ and the structure group
$\Ham\,   (M,\omega)$.
Let us briefly recall the following basic definitions (see
\cite{GLS} for details).

%\smallskip
\begin{defin}
\label{def-Ham-fibr} {\rm A closed 2-form $\tilde{\omega}$ on the
total space $\pi: \Sigma\times M$ of the bundle $\Sigma\times M\to
\Sigma$ is called {\it fiber compatible } if its restriction on
each fiber of $\pi$ is $\omega$. }
\end{defin}
%\bigskip

Let us fix a trivialization of
the bundle $\pi: \Sigma\times M\to \Sigma$ and
let $pr_M :\Sigma\times M\to M$ be the natural projection. The
{\it weak coupling construction} \cite{GLS}, which goes back to
W.Thurston, gives rise to
the following fact.

\begin{prop}[\cite{GLS}]
\label{prop-weak-coupling}
Let $\tilde{\omega}$ be a closed fiber compatible form on $\Sigma\times
M$ that coincides with $pr_M^\ast \omega$ at infinity. Then for a
sufficiently small $\varepsilon > 0$ there exists a smooth family
of closed 2-forms $\{ \Omega_\tau \}$, $\tau \in [0 , \varepsilon]$,
on $\Sigma\times M$ with the following properties:

%\medskip
\noindent (i) $\Omega_0 = \pi^\ast \Omega$;

%\smallskip
\noindent (ii) if one rescales the ends of $\Sigma$
so that $\Sigma$ becomes a compact surface with boundary then
$[\Omega_\tau ] = \tau [\tilde{\omega}] + [\pi^\ast
\Omega]$, where the cohomology classes are taken in
the relative
cohomology group
$H^2
(\Sigma\times M, \partial\Sigma\times M )$;

%\smallskip
\noindent (iii) the restriction of $\Omega_\tau$ on each fiber of
$\pi$ is a multiple of the symplectic form on that fiber;

%\smallskip
\noindent (iv) $\Omega_\tau$ is symplectic for any
$\tau\in (0,\varepsilon]$.
\end{prop}
%\bigskip

%\smallskip
\begin{defin}
\label{def-size} {\rm Define ${\hbox{\it size}}\,
(\tilde{\omega})$ as the supremum of all $\varepsilon$ for which
there exists
a family $\{ \Omega_\tau \}$, $\tau\in [0 , \varepsilon]$,
satisfying the properties (i)-(iv) listed above. }
\end{defin}
%\bigskip

Any fiber compatible closed form $\tilde{\omega}$ defines a
connection $\nabla$ on $\pi: \Sigma\times M\to \Sigma$ whose parallel
transports (with respect to a fixed trivialization) belong to
$\Ham\,  (M,\omega)$. Such a connection is called {\it Hamiltonian}.
Conversely, any Hamiltonian connection $\nabla$ on $\Sigma\times
M\to\Sigma$ can be defined by means of a unique fiber compatible closed
2-form $\tilde{\omega}_\nabla$ such that the 2-form on $\Sigma$
obtained from $\tilde{\omega}_\nabla^{n+1}$ by fiber integration
is 0.
The curvature of a Hamiltonian connection
$\nabla$ can be viewed as a 2-form
associating to each pair $v,w\in T_x\Sigma$ of tangent vectors
a normalized Hamiltonian function $H_{v,w}$ on the fiber
$\pi^{-1} (x)$. The form $\tilde{\omega}_\nabla$ restricted to the
horizontal lifts of vectors $v,w\in T_x\Sigma$ at a point $y\in
\pi^{-1} (x)$ coincides with $H_{v,w} (y)$ \cite{GLS}.

Let $H= (H_1,\ldots , H_l)$ be (time-dependent) {\it
normalized} Hamiltonians on $M$. Let ${\cal C}_H = ({\cal
C}_{H_1}, \ldots , {\cal C}_{H_l})$ be the conjugacy classes in
$\tHam   (M,\omega)$ containing, respectively, the elements
$\widetilde{\varphi}_{H_1},\ldots , \widetilde{\varphi}_{H_l}$.

Fix a trivialization of $\Sigma\times M\to\Sigma$. Let $\Phi_i$,
$i=,1\ldots, l$, and $K>0$ be as in Definition~\ref{def-cal-L}.
Denote by $\widetilde{\cal L} ({\cal C}_H)$ the set of all
Hamiltonian connections $\nabla$ on $\Sigma\times M\to\Sigma$
such that,
with respect to the fixed trivialization, the holonomy of $\nabla$
along the path $\tau \mapsto \Phi_i (s\times e^{2\pi
i\tau})$, $0\leq\tau\leq t$, is $\varphi_{H_i}^t$
for any $s\geq K$, $t\in [0,1]$ and $i=1,\ldots, l$.
Define ${\cal F} ({\cal C}_H )$ as the set of all the forms
$\tilde{\omega}_\nabla$, $\nabla\in \widetilde{\cal L} ({\cal C}_H )$.

%\smallskip
\begin{defin}
\label{def-sigma-weak-coupl} {\rm Define the number
$0< {\hbox{\it size}}_g\, (H) \leq +\infty$
as
\[
{\hbox{\it size}}_g\, (H) = \sup_{\tilde{\omega}\in {\cal F}
({\cal C}_H )} {\hbox{\it size}}\, (\tilde{\omega}).
\]
}
\end{defin}
\bigskip

The following theorem can be proven by exactly the same arguments
as the similar results in \cite{Pol1}, \cite{Pol3} (cf. \cite{ME}).

\begin{prop}
\label{prop-K-area-coupl-Ham} Let ${\hbox{\it K-area}}$ and
$\Upsilon_{l,g}$
be
measured with respect to the bi-invariant Hofer
pse\-udo-metric on
$G=\tHam   (M,\omega)$.
Let $H = (H_1, \ldots , H_l)$ be normalized
Hamiltonians. Then
\[
{\hbox{\it K-area}}_{l,g}\, ({\cal C}_H)\leq {\hbox{\it size}}_g\,
(H),
\]
and, in view of
Proposition~\ref{prop-upsilon-k-area-pos-genus},
\[
\Upsilon_{l,g}\, ({\cal C}_H)\geq 1/{\hbox{\it size}}_g\,
(H).
\]
(If ${\hbox{\it size}}_g\,
(H) = \infty$ we set $1/{\hbox{\it size}}_g\, (H) =0$).
\end{prop}
%\smallskip

%\bigskip
%\bigskip
\section{Pseudo-holo\-morphic curves and an estimate on $\Upsilon_{l,g}$}
\label{sec-upsilon-widetilde-Ham-psh-curves}

In this section we make a crucial step towards the proof
of Theorem~\ref{thm-spectrum-length}:
we obtain an estimate below on $\Upsilon_{l,g}$ for $\tHam  (M,\omega)$ based
on the existence of some pseudo-holomorphic
curves. The relation between existence of such curves and the
hypothesis of Theorem~\ref{thm-spectrum-length} will be discussed
in Section~\ref{sect-pss-stuff}.

{\sl Throughout the section
$(M,\omega)$ can be assumed to be an arbitrary closed connected
symplectic manifold.}

%\bigskip
\subsection{The spaces
${\cal T} (H)$, ${\cal T}_\tau (H)$, ${\cal T}_{\tau, J} (H)$ of
almost complex structures}
\label{subsect-acs-cal-T}

With $\Sigma$ as in the previous section,
consider again the trivial (and  trivialized)
bundle $\pi: \Sigma\times M\to\Sigma$.
Let $j$ be a complex structure on $\Sigma$ compatible with
the area form $\Omega$.
Without loss of
generality we may assume that the identifications $\Phi_i:
[0,+\infty)\times S^1\to \Sigma_i$, $1\leq i\leq l$, are chosen in
such a way that near infinity the complex structure on the ends
gets identified with the standard complex structure on the
cylinder $[0,+\infty) \times S^1$.
Let $J$ be an almost complex structure on $M$ compatible with
$\omega$.

We say that an almost complex structure
$\tilde{J}$
on $\Sigma\times M$ is
$J{\hbox{\it -fibered}}$ if the following conditions are
fulfilled:

\begin{itemize}
\item{$\tilde{J}$ preserves the tangent spaces to
the fibers of $\pi$;}

\item{the restriction of $\tilde{J}$ on any fiber of $\pi$ is
an almost complex structure compatible with the symplectic form
$\omega$ on that fiber;}

\item{the restriction of $\tilde{J}$ to any fiber
$\pi^{-1} (x)$ for $x$ outside of some compact subset of $\Sigma$
is $J$.}
\end{itemize}

Let $H = (H_1,\ldots, H_l)$ be (time-dependent)
Hamiltonians on $M$. Pick ${\hat{\gamma}}\in {\cal
P} (H)$. Let us also pick a cut-off function
$\beta: \R\to [0,1]$ such that $\beta (s)$ vanishes for
$s\leq\epsilon$ and $\beta (s) = 1$ for $s\geq 1-\epsilon$ for
some small $\epsilon>0$.
For each section $u:\Sigma\to \Sigma\times M$ set
\[
u_i := pr_M\circ u\circ \Phi_i : [0,+\infty)\times S^1\to M.
\]
For each $i = 1,\ldots, l$ consider the
non-homogeneous Cauchy-Riemann equation
\begin{equation}
\label{eqn-nonhom-cauchy-riem}
\partial_s u_i + J (u_i)\partial_t u_i - \beta (s) \nabla_u H_i (t,u_i)
= 0,
\end{equation}
where the gradient is taken with respect to the Riemannian metric
$\omega (\cdot , J\cdot )$ on $M$. According to
\cite{Gro-pshc}, the solutions of such an equation
correspond exactly to the pseudo-holomorphic sections of $\pi^{-1}
(\Sigma_i)\to \Sigma_i$ with respect to some unique $J{\hbox{\rm
-fibered}}$ almost complex structure on $\pi^{-1} (\Sigma_i)$.

%\smallskip
\begin{defin}
\label{def-slow-acs-on-a-bundle} {\rm Let $\tilde{J}$ be an almost
complex structure on $\Sigma\times M$ and let $H = (H_1,\ldots ,
H_l)$ be Hamiltonians as above. We shall say that $\tilde{J}$ is
{\it $H{\hbox{\it -compatible}}$ } if there exists  an almost
complex structure $J = J (\tilde{J})$ on $M$ compatible with
$\omega$ such that the following conditions hold:

\begin{itemize}
\item{$\tilde{J}$ is $J{\hbox{\rm -fibered}}$.}

\item{$\pi\circ \tilde{J} = j\circ\pi$.}

\item{There exists a
number $K$ such that for each $i=1\ldots, l$
and each
$\tilde{J}{\hbox{\rm -holomorphic}}$ section
$u$
of $\pi$ the restriction of $u_i : [0, +\infty )
\times S^1\to M$
to $[K,
+\infty )$
is
a solution of the non-homogeneous Cauchy-Riemann equation
(\ref{eqn-nonhom-cauchy-riem}) for $J = J (\tilde{J})$.}
\end{itemize}

Denote by ${\cal T} (H)$ the space of all $H{\hbox{\rm
-compatible}}$ almost complex structures on $\Sigma\times M$. }
\end{defin}
\bigskip

Let ${\cal C}_H = ({\cal C}_{H_1},\ldots , {\cal C}_{H_l})$ be the
conjugacy classes in $\tHam   (M,\omega)$ as in
Section~\ref{sec-Ham-Hofer-metric-def}.

%\bigskip
\noindent
\begin{defin}
\label{def-cal-T-tau} {\rm Consider all the families $\{
\Omega_{\tilde{\omega}_\nabla, \tau} \}$, that arise from the weak
coupling construction associated with $\tilde{\omega}_\nabla$,
$\nabla\in \widetilde{\cal L} ({\cal C}_H)$ (see
Section~\ref{sec-Ham-Hofer-metric-def}). Given a number
$\tau_0\in (0,  {\hbox{\it size}}_g\, (H) )$ consider the set
${\cal Q}_{\tau_0}$ of all the symplectic forms
$\Omega_{\tilde{\omega}_\nabla, \tau_0}$ from the families $\{
\Omega_{\tilde{\omega}_\nabla, \tau} \}$ as above (i.e. we
consider only those families which are defined for the value
$\tau_0$ of the parameter $\tau$ and pick the form
$\Omega_{\tilde{\omega}_\nabla, \tau_0}$ from each such family).
Denote by ${\cal T}_{\tau_0} (H)$ the set of all the almost
complex structures in ${\cal T} (H)$ which are compatible with
some symplec\-tic form from ${\cal Q}_{\tau_0}$.
For an almost complex structure $J$ on $M$ compatible with
$\omega$ denote by ${\cal T}_{\tau, J} (H)$ the set of all
$\tilde{J}\in {\cal T}_\tau (H)$ such that $J = J (\tilde{J})$. }
\end{defin}
%\bigskip

%\bigskip
\subsection{Moduli spaces
${\cal M}_g\, (\hat{\gamma}, H, \tilde{J})$ and the number $s_g
(H)$} \label{subsec-psh-curves-defs}

In \cite{ME} we defined a moduli space ${\cal M} (\hat{\gamma}, H,
\tilde{J})$ of certain pseudo-holomorphic curves of genus
zero. Here we briefly outline the definition modifying it in a
straightforward manner from the case of genus zero to the case
of an arbitrary genus.

In the setup of the previous section
let
\[
\hat{\gamma} =
\biggl[
[\gamma_1, f_1 ],\ldots , [\gamma_l, f_l ]
\biggr]\in {\cal P} (H),
\]
where
$[\gamma_i, f_i ]\in {\cal
P} (H_i)$, $i=1,\ldots ,l$.

Given an almost complex structure $\tilde{J}\in {\cal T} (H)$
consider $\tilde{J}{\hbox{\rm -holomorphic}}$ sections $u:
\Sigma\to \Sigma\times M$ such that the ends
$u_i ([0, +\infty)\times S^1)$, $i=1,\ldots, l$,
of the surface
$pr_M\circ u (\Sigma)\subset M$ converge uniformly at infinity respectively
to the
periodic orbits $\gamma_1,\ldots ,\gamma_l$ and such that
$pr_M\circ u (\Sigma)$ capped off with the discs $f_1 (D^2),
\ldots$, $f_l (D^2)$ is a closed surface representing a torsion
integral homology class in $M$. Denote the space of such
pseudo-holomorphic sections $u$ by ${\cal M}_g\,
(\hat{\gamma}, H, \tilde{J})$ (the index $g$ indicates the genus
of the $\tilde{J}{\hbox{\rm -holomorphic}}$ curves we
consider).

%\bigskip
\noindent
\begin{defin}[cf. \cite{ME}]
\label{def-spec-g-s-g} {\rm We will say that a real number $c$
is $g{\hbox{\it -durable}}$ if there exists
\begin{itemize}
  \item{a sequence
$\{ \tau_k \}\nearrow {\hbox{\it size}}_g\, (H)$,}
  \item{a sequence
$\{ \hat{\gamma}_k\}$, $\hat{\gamma}_k\in {\cal P} (H)$, such that
$\lim_{k\to +\infty} {\cal A}_H (\hat{\gamma}_k ) = c$, }
  \item{a sequence
$\{ \tilde{J}_{\tau_k}\}, \tilde{J}_{\tau_k}\in {\cal T}_{\tau_k}
(H)$,}
\end{itemize}
so that for all $k$ the spaces ${\cal M}_g\, (\hat{\gamma}_k, H,
\tilde{J}_{\tau_k})$ are non-empty.

Denote by $s_g (H)$ the supremum
of all $g{\hbox{\rm -durable}}$ numbers.
If  ${\hbox{\it spec}}_g\, (H)$ is empty let $s_g (H) = -\infty$.
}
\end{defin}

%\bigskip
\subsection{Estimating
$\Upsilon_{l,g}$ by actions of periodic orbits}
\label{subsect-upsilon-l-g-main-estimate}

Let $(M,\omega)$ be a closed connected symplectic manifold.
Consider the function $\Upsilon_{l,g}$ on conjugacy classes in
$\tHam  (M,\omega)$ defined by means of
the Hofer pseudo-metric on
$\tHam  (M,\omega)$
(see
Section~\ref{sec-Ham-Hofer-metric-def}).

%\smallskip
\begin{prop}
\label{prop-exist-psh-curve-implies-estimate-pos-genus} Suppose
the Hamiltonians $H = (H_1,\ldots, H_l)$ are normalized. Then
$\Upsilon_{l,g}\, ({\cal C}_H ) \geq s_g (H)$.
\end{prop}
%\smallskip

For the proof of
Proposition~\ref{prop-exist-psh-curve-implies-estimate-pos-genus}
see
Section~\ref{subsect-pf-psh-curve-leads-to-estimate}.
Proposition~\ref{prop-exist-psh-curve-implies-estimate-pos-genus}
is a generalization of Theorem 1.3.1 from \cite{ME} to the case of
pseudo-holomorphic curves of an arbitrary genus.

%\bigskip
%\bigskip
\subsection{The proof of
Proposition~\ref{prop-exist-psh-curve-implies-estimate-pos-genus}}
\label{subsect-pf-psh-curve-leads-to-estimate}

The proof virtually repeats the proof of Theorem 1.3.1 from
\cite{ME} in the case of genus zero.

Without loss of generality we may assume that $s_g (H) >0$
(otherwise the proposition is trivial).
Fix an arbitrary small $\epsilon > 0$ so that
$s_g (H) -\epsilon >0$.
Then there exists a $g{\hbox{\rm -durable}}$ number $c$ such
that
\[
s_g (H)\geq c\geq s_g (H) -\epsilon >0.
\]
Now, according to Definition~\ref{def-spec-g-s-g}, for any
sufficiently small
$\delta >0$
there exist
\begin{itemize}
\item{a number $\tau_0$, ${\hbox{\it size}}_g\, (H) -\delta \leq
\tau_0 < {\hbox{\it size}}_g\, (H)$},

\item{a Hamiltonian connection
$\nabla\in \widetilde{\cal L} ({\cal C}_H)$ and the corresponding
2-form $\tilde{\omega}_\nabla$ on $\Sigma\times M$ such that
$\tau_0 \leq {\hbox{\it
size}}\, (\tilde{\omega}_\nabla )$,}

\item{a symplectic form
$\{\Omega_{\tilde{\omega}_\nabla , \tau_0 }\}$
from a weak coupling deformation,}

\item{an almost complex structure
$\tilde{J}_{\tau_0}\in {\cal T}_{\tau_0} (H)$
compatible with $\{\Omega_{\tilde{\omega}_\nabla , \tau_0 }\}$,}

\item{$\hat{\gamma}\in {\cal P} (H)$ such that ${\cal A}_H
(\hat{\gamma})\geq c -\delta >0$,}

\end{itemize}
so that the space
${\cal M}_g (\hat{\gamma}, H, \tilde{J}_{\tau_0} )$ is non-empty.

Pick a map $u\in {\cal M}_g (\hat{\gamma}, H, \tilde{J}_{\tau_0} )$.
Then
\[
0\leq \int_{u(\Sigma)}
\Omega_{\tilde{\omega}_\nabla ,\tau_0 }
\]
because
$\tilde{J}_{\tau_0}$ is compatible with the symplectic form
$\Omega_{\tilde{\omega}_\nabla,\tau_0 }$ and the surface
$u (\Sigma)\subset \Sigma\times M$ is a
pseudo-holomorphic curve with respect to
$\tilde{J}_{\tau_0}$.
On the other hand,
\[
\int_{u(\Sigma)}
\Omega_{\tilde{\omega}_\nabla ,\tau_0 } = \int_{u(\Sigma)} \tau_0
\tilde{\omega}_\nabla + \int_{u(\Sigma)}\pi^\ast \Omega
\]
because of the cohomological condition
satisfied by a weak coupling deformation (see
condition (iii) in Proposition~\ref{prop-weak-coupling}).
Thus
\begin{equation}
\label{eqn-ineq-tot-integral-positive} 0\leq \int_{u(\Sigma)} \tau_0
\tilde{\omega}_\nabla + \int_{u(\Sigma)}\pi^\ast \Omega.
\end{equation}

Next we recall the following lemma (see Lemma 5.0.1 in \cite{ME}) whose
proof does not depend on the genus of $\Sigma$.

%\smallskip
\begin{lem}
\label{lem-integral-equal-action}
\begin{equation}
\label{equal-action-omega-tilde} \int_{u(\Sigma)}
\tilde{\omega}_\nabla = - {\cal A}_H (\hat{\gamma}).
\end{equation}
\end{lem}
%\smallskip

Now, using Lemma~\ref{lem-integral-equal-action} and the
fact that the total $\Omega{\hbox{\rm -area}}$ of $\Sigma$ is 1,
one can rewrite (\ref{eqn-ineq-tot-integral-positive}) as
\[
\tau_0\leq \frac{1}{{\cal A}_H (\hat{\gamma})}
\]
and hence
\[
{\hbox{\it size}}_g\, (H) -\delta \leq
\tau_0\leq \frac{1}{{\cal A}_H (\hat{\gamma})}\leq  \frac{1}{c-\delta}.
\]
Since this is true for any $\delta >0$,
\[
{\hbox{\it size}}_g\, (H) \leq \frac{1}{c}\leq \frac{1}{s_g (H)
-\epsilon}.
\]
In view of
Proposition~\ref{prop-K-area-coupl-Ham},
\[
\Upsilon_{l,g} ({\cal C}_H )\geq
s_g (H) -\epsilon.
\]
Since $\epsilon > 0$ was chosen arbitrarily,
\[
\Upsilon_{l,g} ({\cal C}_H )\geq s_g (H)
\]
and the proposition is proven.
\b

%\bigskip
%\bigskip
\section{Pair-of-pants product on Floer cohomology and the moduli spaces
${\cal M}_g\, (\hat{\gamma}, H, \tilde{J})$}
\label{sect-pss-stuff}

The goal of this section is to relate the moduli spaces ${\cal
M}_g\, (\hat{\gamma}, H, \tilde{J})$ and the number $s_g (H)$ to
the hypothesis of Theorem~\ref{thm-spectrum-length}. This will be
done as in \cite{ME} by means of the multiplicative structure on
Floer and quantum cohomology which will be used to guarantee that
certain moduli spaces ${\cal M}_g\, (\hat{\gamma}, H, \tilde{J})$
with ${\cal A}_H (\hat{\gamma}) >0$
are non-empty and for that reason $s_g (H) > 0$.

{\sl Throughout this section
the symplectic manifold $(M^{2n},\omega)$ is assumed to be strongly semi-positive.}

Let $\Sigma$ be a Riemann surface of genus $g$ as above with $l=1$
cylindrical end. Let $(H, J)$ be a regular Floer pair. Fix a number
$\tau$, $0 < \tau < {\hbox{\it size}}_g\, (H)$.

%\smallskip
\begin{prop}
\label{prop-moduli-space-generically-smooth} For a gene\-ric
$\tilde{J}_\tau \in {\cal T}_{\tau ,J} (H)$ and any
$\hat{\gamma}\in {\cal P} (H)$ with the Conley-Zehnder index $\mu
(\hat{\gamma}) = 2n (1-g)$ the space ${\cal M}_g (\hat{\gamma}, H,
\tilde{J}_\tau )$ is either empty or an oriented compact
zero-dimensional manifold.
\end{prop}
%\smallskip

Given such a generic $\tilde{J}_\tau$ we will say that the pair
$(H,\tilde{J}_\tau )$ is {\it regular. } For a regular pair $(H,
\tilde{J}_\tau )$, $\tilde{J}_\tau \in {\cal T}_{\tau ,J} (H)$,
and an element $\hat{\gamma}\in {\cal P} (H)$,
$\mu({\hat{\gamma}}) = 2n (1-g)$, count the curves from the
compact oriented zero-dimensional moduli space ${\cal M}_g
(\hat{\gamma}, H, \tilde{J}_\tau )$ with their signs. The
resulting Gromov-Witten number will be denoted by $\# {\cal M}_g (\hat{\gamma}, H,\tilde{J}_\tau )$. Take the sum
\begin{equation}
\label{eqn-Floer-cycle} \sum_{\hat{\gamma}} \# {\cal M}_g (\hat{\gamma}, H,\tilde{J}_\tau ) \hat{\gamma},
\end{equation}
over all $\hat{\gamma}\in {\cal P} (H)$ such that $\mu
(\hat{\gamma}) = 2n (1-g)$. The sum in (\ref{eqn-Floer-cycle})
represents an integral chain $\theta_{\Sigma, H,\tilde{J}_\tau }$
in the chain complex $CF_\ast (H, J)$ or, from the
Poincar{\'e}-dual point of view, an integral cochain
$\theta^{\Sigma, H,\tilde{J}_\tau }$ in the cochain complex
$CF^\ast (\bar{H}, J)$.

The following proposition from \cite{ME} is a minor generalization
of Theorem 3.1 in \cite{PSS}: we use a bigger class of admissible
almost complex structures but the proof can be carried out in
exactly the same way (see \cite{ME} for a discussion
and an outline of the proof).

%\smallskip
\begin{prop}
\label{prop-from-main-pss} Let $(H, J)$ be a regular Floer pair.
For any $\tau$, $0< \tau <
{\hbox{\it size}}_g\, (H)$, and a regular pair $(H, \tilde{J}_\tau
)$ as above the cochain $\theta^{\Sigma, H,\tilde{J}_\tau }$
defines a cocycle in the cochain complex $CF^\ast\, (\bar{H}, J)$.
The corresponding cohomology class
\[
\Theta^{H, J, g}\in HF^\ast (\bar{H}, J)
\]
is of degree $2ng$ and does not depend on $\tau$ and
$\tilde{J}_\tau \in {\cal T}_{\tau, J} (H)$. Likewise the chains
$\theta_{\Sigma, H,\tilde{J}_\tau }$ in $CF_\ast (H, J)$ are
cycles and represent a homology class
\[
\Theta_{H, J, g}\in HF_\ast (H, J),
\]
of degree $2n(1-g)$ which does not depend on $\tau$ and
$\tilde{J}_\tau \in {\cal T}_{\tau, J} (H)$ and is
Poincar{\'e}-dual to $\Theta^{H, J, g}$.

\end{prop}
%\smallskip

%\smallskip
\begin{prop}
\label{prop-Theta-represents-fundam-class}
Let $E\in QH^{2n} (M,\omega)$ be the Euler class and
let $(H, J)$ be a regular Floer pair.
Then
\[
\Psi_{{\bar H}, J}^{-1} (\Theta^{H, J, g} ) = E^g \in QH^{2ng}
(M,\omega).
\]
\end{prop}
%\smallskip

Postponing the proof of
Proposition~\ref{prop-Theta-represents-fundam-class} we first
state the main result of this section.

%\bigskip
\noindent
\begin{prop}
\label{prop-s-g-geq-c-e-g-minus-omega-A-g} Suppose that
$H$ belongs to a regular Floer pair $(H, J)$ and that the class
$E^g \in QH^{2ng} (M,\omega)$ is non-zero. Then
$s_g (H)\geq c (E^g, H)$.
\end{prop}
%\smallskip

%\bigskip
\noindent {\bf Proof of
Proposition~\ref{prop-s-g-geq-c-e-g-minus-omega-A-g}.}

\noindent
According to
Proposition~\ref{prop-Theta-represents-fundam-class},
\[
PD (\Psi_{\bar{H}, J} (E^g )) = \Theta_{H, J, g}
\]
and, according to Proposition~\ref{prop-from-main-pss},
for any $\tau$, $0 < \tau < {\hbox{\it
size}}_g\, (H)$, the Floer homology class $\Theta_{H, J, g}$ can
be represented by a chain
\begin{equation}
\label{eqn-Floer-cycle-1} \theta_{\Sigma, H,\tilde{J}_\tau } =
\sum_{\hat{\gamma}} \# {\cal M}_g (\hat{\gamma}, H,\tilde{J}_\tau )
\hat{\gamma},
\end{equation}
from $CF_\ast (H, J)$, where the sum in (\ref{eqn-Floer-cycle-1})
is taken over all $\hat{\gamma}\in {\cal P} (H)$ such that $\mu
(\hat{\gamma}) = 2n (1-g)$. According to part 4 of
Proposition~\ref{prop-properties-c-general-case},
for any $\epsilon
> 0$ the sum (\ref{eqn-Floer-cycle-1}) has to contain a non-zero term $\# {\cal M}_g (\hat{\gamma}, H,\tilde{J}_\tau ) \hat{\gamma}$ such that ${\cal A}_H
(\hat{\gamma}) > c (E^g , H) -\epsilon$. Recalling
Definition~\ref{def-spec-g-s-g} we see that the number
$c (E^g, H)$ is $g{\hbox{\rm -durable}}$ and therefore
$s_g (H)$, which is the supremum of all
$g{\hbox{\rm -durable}}$ numbers, is no less than
$c (E^g , H)$.
The proposition is proven.\ \b

\bigskip
\noindent
{\bf Proof of
Proposition~\ref{prop-Theta-represents-fundam-class}}
%\label{subsect-pf-prop-Theta-represents-fundam-class}

\noindent
For brevity denote $V_k = QH_k (M, \omega), \bar{V}^k = QH^k (M,
\omega) = {\hbox{\it Hom}} (V_k, \Lambda_\omega^0)$, $V = \oplus
V_k, \bar{V} = \oplus \bar{V}^k$. Without loss of generality we
can identify the Floer cohomology with $\bar{V}$ and the Floer
homology with $V$ by means of the Piunikhin-Salamon-Schwarz
isomorphisms as above. Consider the spaces of the form $V^{\otimes l_1}
\otimes {\bar{V}}^{\otimes l_2}$ where $\otimes$ stands for tensor
product over $\Lambda_\omega^0$.

The pair-of-pants product on the Floer cohomology
can be viewed as a cohomological operation on $HF^\ast (H, J)$
associated to a surface of genus zero with two "entering"
and one "exiting" cylindrical ends.
In fact it is a
part of a more general series of cohomological operations on
$HF^\ast (H, J)$ defined by means of the moduli spaces ${\cal M}_g\,
(\hat{\gamma}, H, \tilde{J})$ for surfaces of an arbitrary genus
$g$ and with an arbitrary number of cylindrical ends. This was first
proven \cite{Sch-PhD} in a simpler setup but can be shown in our
case in exactly the same way (cf. \cite{PSS}, also see \cite{ME}).

More precisely,
every Riemann surface of a
genus $g\geq 0$ with $l_1\geq 0$ positively oriented
ends ("entrances") and $l_2\geq 0$ negatively oriented ends
("exits") gives rise to a certain element
\[
\Xi_{l_1, l_2, g}\in V^{\otimes l_1} \otimes {\bar{V}}^{\otimes
l_2},
\]
which depends only on $l_1, l_2, g$. In case of a closed surface,
when $l_1 = l_2 = 0$, the element $\Xi_{0, 0, g}\in
\Lambda_\omega^0$ is an {\it integer number} -- it comes
as a result of counting certain {\it closed} pseudo-holomorphic
curves with {\it integral multiplicities} (see \cite{PSS},
\cite{Sch-PhD}).

Each element $\Xi_{l_1, l_2, g}$ can be viewed as a polylinear map
over $\Lambda_\omega^0$, denoted by $\xi_{l_1, l_2, g}$, which
sends a tuple of $l_1$ quantum cohomology classes $(f_1,\ldots,
f_{l_1})$ from $\bar{V}$ to an element of ${\bar{V}}^{\otimes
l_2}$. Namely, each monomial term
\[
\alpha_1\otimes\ldots\otimes\alpha_{l_1}\otimes
g_1\otimes\ldots\otimes g_{l_2},
\]
$\alpha_i\in V$, $g_i\in \bar{V}$, in $\Xi_{l_1, l_2, g}$ sends
$f_1\otimes\ldots \otimes f_{l_1}\in {\bar{V}}^{\otimes l_1} $ to
\[
\prod_{i=1}^{l_1} (f_i, \alpha_i) \cdot  g_1\otimes\ldots\otimes
g_{l_2},
\]
where $(\cdot, \cdot)$ is the evaluation pairing with values in
$\Lambda_\omega^0$ (see Section~\ref{subsec-Novikov-quantum}). In
particular, if $f_i\in V^k$ and $\alpha_i\in V_j$ then $(f_i,
\alpha_i) = 0$ unless $k=j$.

We will use the Poincar{\'e} isomorphism $PD: V\to \bar{V}$ to
raise and lower indices of tensors:
\[
PD: V^{\otimes l_1} \otimes {\bar{V}}^{\otimes l_2} \to V^{\otimes
\{l_1 -1\} } \otimes {\bar{V}}^{\otimes \{ l_2 +1 \} }.
\]
Here the $l_1$-th factor $V$ in the product $V^{\otimes l_1}
\otimes {\bar{V}}^{\otimes l_2}$ is transformed by $PD$ into
$\bar{V}$.

The basic "topological field theory" properties of $\xi_{l_1, l_2,
g}, \Xi_{l_1, l_2, g}$ and the relation between them and the ring
structure of the quantum cohomology can be described by the
following proposition:

%\bigskip
\begin{prop}[\cite{PSS},\cite{Sch-PhD}]
\label{prop-tft-properties} The following properties hold for any
$l_1$, $l_2$, $l_3$, $g, g_1, g_2$:

\medskip
\noindent 1) $\xi_{l_1, l_2, g_1} \circ \xi_{l_2, l_3, g_2} =
\xi_{l_1, l_3, g_1 + g_2}$ ("topological field theory property").

\medskip
\noindent 2) $PD (\Xi_{l_1, l_2, g}) = \Xi_{l_1 -1, l_2+1, g}$, if
$l_1\geq 1$.

\medskip
\noindent 3) $\Xi_{0,1,0} = {\bf 1} \in \bar{V}^0$ and
$\Xi_{1,0,0}$ is the class $PD ({\bf 1})\in V_{2n}$,
Poincar{\'e}-dual to ${\bf 1}$ in $V$. Thus the map
$\xi_{1,0,0}:\bar{V} \to \Lambda_\omega^0$ is the evaluation
pairing with $PD ({\bf 1})\in V_{2n}$: its sends each $\beta =
\sum_i \beta_i\in \bar{V}$, $\beta_i\in \bar{V}^i$, to
$\xi_{1,0,0} (\beta ) = (\beta_{2n}, PD ({\bf 1}))$.

\medskip
\noindent 4) $\xi_{1,1,0} = Id: \bar{V}\to \bar{V}$.

\medskip
\noindent 5) $\xi_{2,1,0}: \bar{V} \otimes \bar{V} \to
\bar{V}$ is the quantum (pair-of-pants) multiplication.

\end{prop}
%\smallskip

The class $\Psi_{\bar{H}, J}^{-1} (\Theta^{H, J, g} )\in QH^\ast (M,
\omega)$ is precisely the element $\Xi_{0,1,g}\in \bar{V}$
corresponding to a surface of genus $g$ with one "exiting" end
(see \cite{PSS}, \cite{Sch-PhD}). On the other hand, such a surface
can be obtained by taking a surface
with $g$ "entrances" and one "exit" and gluing to each "entrance" a
surface of genus 1 with one "exit". In view of parts 1 and 5 of
Proposition~\ref{prop-tft-properties} we only need
to prove the following lemma.

%\bigskip
\begin{lem}
\label{lem-Xi-0-1-1} $\Xi_{0,1,1} = E.$
\end{lem}

%\bigskip
\noindent {\bf Proof of Lemma~\ref{lem-Xi-0-1-1}.}

\noindent We will present the cohomology class $\Xi_{0,1,1}$,
viewed as a $\Lambda_\omega^0{\hbox{\rm -linear}}$ map $V\to
\Lambda_\omega^0$, as a composition of two polylinear forms {\it
over $\Lambda_\omega$.} This will allow us to compute
$\Xi_{0,1,1}$ by means of $\Lambda_\omega{\hbox{\rm -bases}}$ of
$\bar{V}$.

Let $\{ e_i\}$, $e_i\in \bar{V}^{k_i}$, be a
$\Lambda_\omega{\hbox{\rm -basis}}$ of $\bar{V}$ and let $\{
F_i\}$, $F_i\in V_{k_i}$, be the dual $\Lambda_\omega{\hbox{\rm
-basis}}$ of $V$ so that $(e_i, F_j) =\delta_{ij}$. Set $\bar{e}_i
: = PD (F_i)\in \bar{V}^{2n-k_i}$. Then $\{
\bar{e}_i\}$ is a $\Lambda_\omega{\hbox{\rm -basis}}$ of $\bar{V}$
Poincar{\'e}-dual to the basis $\{ e_i\}$ (see
Section~\ref{subsec-Novikov-quantum}).

Now recall part 4 of Proposition~\ref{prop-tft-properties} and
observe that the map $\xi_{1,1,0} = Id: \bar{V}\to \bar{V}$ is
obviously $\Lambda_\omega{\hbox{\rm -linear}}$. Thus $\Xi_{1,1,0}$
is, in fact, an element of $V\otimes_{\Lambda_\omega} \bar{V}$,
where $\otimes_{\Lambda_\omega}$ stands for the graded tensor product {\it
over $\Lambda_\omega$} as opposed to $\otimes$ which denotes
the tensor product over $\Lambda_\omega^0$.
(Recall that as far as the grading is concerned, $\bar{V}$ and $V$ are
graded modules over, respectively, the
graded and the anti-graded versions of $\Lambda_\omega$ --
see
Section~\ref{subsec-Novikov-quantum}).
Hence, using the chosen bases of $V$ and
$\bar{V}$ over $\Lambda_\omega$, one can write
\[
\Xi_{1,1,0} = Id =\sum_i F_i\otimes e_i\in V\otimes_{\Lambda_\omega} \bar{V}.
\]

Since the Poincar{\'e} isomorphism $PD$ is $\Lambda_\omega{\hbox{\rm
-linear}}$, the fact that $\Xi_{1,1,0} \in V\otimes_{\Lambda_\omega} \bar{V}$
together with part 2 of Proposition~\ref{prop-tft-properties}
tells us that $\Xi_{0,2,0}\in \bar{V}\otimes_{\Lambda_\omega} \bar{V}$. With
respect to the chosen basis of $\bar{V}$ over $\Lambda_\omega$ it
can be written as
\begin{equation}
\label{eqn-xi-0-2-0} \Xi_{0,2,0} = \sum_i \bar{e}_i \otimes e_i\in
\bar{V}\otimes_{\Lambda_\omega} \bar{V}.
\end{equation}

The quantum multiplication $\xi_{2,1,0}: \bar{V}\times \bar{V}\to
\bar{V}$ is also $\Lambda_\omega{\hbox{\rm -linear}}$ because it
is associative \cite{Liu}, \cite{McD-Sal-pshc}, \cite{Ru-Ti},
\cite{Ru-Ti-1}. Thus $\Xi_{1,1,0}\in
V\otimes_{\Lambda_\omega} V\otimes_{\Lambda_\omega} \bar{V}$.

Part 1 of Proposition~\ref{prop-tft-properties} tells us that
$\Xi_{0,1,1}$ can be found if we compose the operation associated
to a surface of genus 0 with two "exits" and the operation
associated with a surface of genus 0 with 2 "entrances" and 1
"exit". Since, as we have explained above, both operations are, in
fact, $\Lambda_\omega{\hbox{\rm -linear}}$, we can use
(\ref{eqn-xi-0-2-0}) together with parts 1 and 5 of
Proposition~\ref{prop-tft-properties} and express $\Xi_{0,1,1}\in
\bar{V}$ in terms of the chosen basis of $\bar{V}$ {\it over
$\Lambda_\omega$}:
\[
\Xi_{0,1,1} = \sum_i \bar{e}_i \ast e_i = \sum_i (-1)^{{\rm deg}\,
e_i} e_i \ast \bar{e}_i = E\in \bar{V}.
\]
This finishes the proof of
Lemma~\ref{lem-Xi-0-1-1} and
Proposition~\ref{prop-Theta-represents-fundam-class} are proven.
\b

%\smallskip
\begin{rem}
\label{rem-schwarz-euler-char} {\rm Recall from
Section~\ref{subsec-Novikov-quantum}
that the component of degree $2n$ of $E$ is equal to $\chi (M) m\in
H^{2n} (M, \Q)$. Thus, applying parts 1 and 3 of
Proposition~\ref{prop-tft-properties} one gets that $\Xi_{0, 0, 1}
= \chi (M)$ -- see Corollary  5.4.12 from \cite{Sch-PhD}. }
\end{rem}
%\smallskip

%\bigskip
%\bigskip
\section{Proof of
Theorem~\ref{thm-spectrum-length}}
\label{sect-pf-thm-spectrum-length-part-II}

If $H$ belongs to a regular Floer pair then
Proposition~\ref{prop-comm-length-upsilon}
(in the case $l=1$) and
Proposition~\ref{prop-s-g-geq-c-e-g-minus-omega-A-g}
yield
\[
\Upsilon_{1,g} ({\cal C}_H)\geq c (E^g, H).
\]
Since both sides depend continuously on $H$, the inequality
is in fact true for an arbitrary $H$. Thus if $c (E^g, H)$
then $\Upsilon_{1,g} ({\cal C}_H) >0$.
In view of Proposition~\ref{prop-comm-length-upsilon}
and
(\ref{eqn-ham-symp-2}), it shows that ${\hbox{\it
cl}}\,(\widetilde{\varphi}_H) >g$ and the theorem is proven.

For the readers benefit we now quickly review the course of the
proof. We used the function $\Upsilon_{1,g}$ on conjugacy classes
in
$\tHam  (M,\omega)$ to measure the distance, with respect to the Hofer pseudo-metric,
from the conjugacy class ${\cal C}_H$
of an element $\widetilde{\varphi}_H\in \tHam  (M,\omega)$
to the set
of elements whose commutator length does not exceed $g$.
We passed to the description of $\Upsilon_{1,g}$ in terms of K-area
and used the weak coupling construction
to estimate $\Upsilon_{1,g} ({\cal C}_H)$
from below by a non-negative (but
possibly zero) number
$1/{\hbox{\it size}}_g\, (H)$.
Thus, as soon as $1/{\hbox{\it size}}_g\, (H)$
is non-zero,
the commutator length of $\widetilde{\varphi}_H$ is greater than
$g$.

In order to guarantee that
${\hbox{\it size}}_g\, (H) < +\infty$
and
$1/{\hbox{\it size}}_g\, (H) > 0$
we first showed that any number
$0 <\tau <{\hbox{\it size}}_g\, (H)$
can be estimated from above by
$1/{\cal A}_H (\hat{\gamma})$
if for a connected oriented surface $\Sigma$
of genus $g$ with one cylindrical end the moduli space
${\cal M} (\hat{\gamma}, H,
\tilde{J}_\tau)$ is non-empty for a certain almost complex
structure  $\tilde{J}_\tau$ compatible with a symplectic form
corresponding to the parameter $\tau$ in some weak coupling
deformation.

Then we discussed how the condition
$c (E^g, H) > 0$
would guarantee the existence of elements $\hat{\gamma}$
as above and provide
a uniform positive bound from below on their actions
as
$\tau$ tends to ${\hbox{\it size}}_g\, (H)$ -- such a positive
bound would also estimate the number
$1/{\hbox{\it size}}_g\, (H)$ from below and
makes sure it is non-zero.
In order to deduce the existence of the bound from the condition
$c (E^g, H) > 0$, assume without loss of generality that
$H$
belongs to a regular Floer pair.
Pick a generic
$\tilde{J}_\tau$
so that all
the moduli space
${\cal M}_g (\hat{\gamma}, H,
\tilde{J}_\tau)$
of expected dimension zero are either empty or smooth, compact and
oriented.
Consider all the elements $\hat{\gamma}\in {\cal P} (H)$
which give rise to such
${\cal M}_g (\hat{\gamma}, H,
\tilde{J}_\tau)$. The
sum
$\sum \# {\cal M}_g (\hat{\gamma}, H,
\tilde{J}_\tau) \hat{\gamma}$
over all such $\hat{\gamma}$ is
a Floer cycle which represents a Floer homology class
corresponding, under the composition of
Poincar{\'e} and Piunikhin-Salamon-Schwarz isomorphisms, to
the quantum cohomology class $E^g$. The condition
$c (E^g, H) >0$ guarantees that
the sum
$\sum \# {\cal M}_g (\hat{\gamma}, H,
\tilde{J}_\tau) \hat{\gamma}$
contains at least one non-zero term corresponding to some
$\hat{\gamma}$ with a positive action arbitrarily close to
$c (E^g, H)$ from below. Since $\tau$ can be chosen arbitrarily
close to ${\hbox{\it size}}_g\, (H)$
this construction provides the necessary
positive estimate from below on $1/{\hbox{\it size}}_g\, (H)$ --
as long as $c (E^g, H) > 0$.

Thus if $c (E^g, H) > 0$ the distance
$\Upsilon_{1,g} ({\cal C}_H)$  is
positive and therefore
the commutator length of $\widetilde{\varphi}_H$ is greater than
$g$.\ \b

%\bigskip
%\bigskip
\section{Proof of
Theorem~\ref{thm-comm-length-closed-mfds}}
\label{sect-quasimorphisms}

As in \cite{Ba-Ghys} the construction of ${\mathfrak f}$ will
involve two main ingredients: a homogeneous quasimor\-phism
$\tau$ on the universal cover $\widetilde{\hbox{\it Sp}}\, (2n,
\R)$ of the group ${\hbox{\it Sp}}\, (2n, \R)$ and a
collection of functions $F_{\{\phi_t\}}: M\to \widetilde{\hbox{\it
Sp}}\, (2n, \R)$ that will be defined for any path $\{\phi_t
\}$ in $\Symp_0\, (M,\omega)$. The only place where one needs
to make adjustments to the argument from \cite{Ba-Ghys} is the
construction of $F_{\{\phi_t\}}$. The quasimorphism
$\tau$ that we will use is the same as in \cite{Ba-Ghys} -- we
recall its definition below.

%\bigskip
\subsection{The quasimorphism $\tau$
on $\widetilde{\hbox{\it Sp}}\, (2n, \R)$}
\label{subsect-lagr-grassm}

Let $\R^{2n}$ be the standard linear symplectic space with the symplectic
form
$dp\wedge dq$ on it.
Let
$\Lambda (\R^{2n})$ be the {\it Lagrange Grassmannian} of
$\R^{2n}$, i.e. the space of all Lagrangian planes in $\R^{2n}$.
It is a
compact manifold that can be identified with $ U (n)/ O (n)$.
The map associating to a unitary
matrix the square of its  determinant descends to
a map: $det^2 : \Lambda (\R^{2n})\to S^1$.
Set the Lagrangian $p{\hbox{\rm -coordinate}}$ plane $L_0\subset \R^{2n}$
as a base point in $\Lambda (\R^{2n})$.
Fix a lift
$\widetilde{L}_0$ of $L_0$ in the
universal cover $\widetilde{\Lambda} (\R^{2n})$ of $\Lambda
(\R^{2n})$ corresponding to the constant path in $\Lambda
(\R^{2n})$ identically equal to $L_0$. Let
$\widetilde{det}^2 : \widetilde{\Lambda} (\R^{2n})\to \R$ be the lift of $det^2$ such that $\widetilde{det}^2
(\widetilde{L}_0 ) = 0$. Roughly speaking, to any
path of Lagrangian subspaces starting at $L_0$ the
function $\widetilde{det}^2$ associates its rotation number in
$\Lambda (\R^{2n})$ (with the whole construction depending on
our choice of the Darboux basis in $\R^{2n}$).

The group $Sp\, (2n, \R)$ acts transitively on $\Lambda (\R^{2n})$ and the universal cover $\widetilde{Sp}\, (2n, \R)$
acts on $\widetilde{\Lambda} (\R^{2n})$. The fundamental
group of $\Lambda (\R^{2n})$ is isomorphic to ${\bf Z}$ and
the map induced by the projection $U (n) \to U (n)/ O (n) =
\Lambda (\R^{2n})$ on $\pi_1 (U (n))$ is the multiplication
by two. On the other hand $U (n)$ is a deformation retract of
$Sp\, (2n, \R)$. These identifications of $\pi_1
(\widetilde{\Lambda} (\R^{2n}) )$ and
$\pi_1 (Sp\, (2n, \R))$ with ${\bf Z}$ by means of
$\pi_1 (U (n))$ lead
to the following proposition.

%\smallskip
\begin{prop}
\label{prop-boundedness-sp-lagr-grassm} A subset of
$\widetilde{Sp}\, (2n, \R)$ is bounded if and only if the
function
\[
\widetilde{\Phi}\mapsto \widetilde{det}^2 (\widetilde{\Phi}
(\widetilde{L}_0))
\]
is bounded on it.
\end{prop}
%\smallskip

Now for $\widetilde{\Phi}\in \widetilde{Sp}\, (2n, \R)$ set
\[
\tau_{det} (\widetilde{\Phi} ) := \widetilde{det}^2
(\widetilde{\Phi} (L_0))
\]
and define
\[
\tau (\widetilde{\Phi}) := \lim_{k\to +\infty} \frac{\tau_{det}
(\widetilde{\Phi}^k)}{k}.
\]

%\smallskip
\begin{prop}[\cite{Ba-Ghys}]
\label{prop-tau-barge-ghys} The function $\tau : \widetilde{Sp}\,
(2n, \R)\to \R$ is continuous and, moreover, it is a homogeneous
quasimorphism on $\widetilde{Sp}\, (2n, \R)$.

\end{prop}

%\bigskip
\subsection{Construction and basic properties of
$F_{\{\phi_t\}}$}

%\smallskip
Let $\mu$ be the measure on $M^{2n}$ defined by the volume form
$\omega^n$. Pick an almost complex structure $J$ on $M$ compatible
with $\omega$.
The Riemannian metric $\omega (\cdot, J \cdot)$ and the form $\omega$ can be viewed as the
real and the imaginary part of a Hermitian metric on $M$.

%\smallskip
\begin{lem}
\label{lem-fuchs} Let $(M^{2n},\omega)$ be a closed symplectic manifold
with $c_1 (M) = 0$. Then there exists a compact triangulated
subset $Y\subset M$ of codimension at least 3
such that the tangent bundle of $M$ admits a unitary
trivialization over $M\setminus Y$.
\end{lem}

%\bigskip
\noindent {\bf Proof of Lemma~\ref{lem-fuchs}.}

\noindent Start with a piecewise smooth triangulation of $M$. Since $M$
is orientable and $c_1 (M) = 0$ there exists a unitary
trivialization of $T M$ over the 2-skeleton of $M$. Consider
the barycentric star decomposition associated with the
triangulation. Let
$Y$ be the $(2n-3){\hbox{\rm -skeleton}}$ of the barycentric star
cell decomposition. Then the 2-skeleton of the original triangulation
of $M$ is a
deformation retract of $M\setminus Y$. Therefore there exists a
unitary trivialization of $T M$ over $M\setminus Y$. The
lemma is proven.
\b
%\medskip

%\smallskip
\begin{rem}
\label{rem-lem-fuchs} {\rm
Obviously the unitary trivialization from
Lemma~\ref{lem-fuchs} is also symplectic: it symplectically
identifies each tangent space of $M\setminus Y$ with the standard
symplectic space $\R^{2n}$.
}
\end{rem}
%\smallskip

Let ${\cal F} \to M$ be the bundle of all
unitary frames of $T M$.
Use the Riemannian metric on $M$ and a metric on $U (n)$
to define a metric on the total space ${\cal F}$.
Any unitary trivialization of
$TM$ over $M\setminus Y$ can be viewed
as a section
$f$
of ${\cal F} \to M$ over $M\setminus Y$. The Riemannian
metrics on ${\cal F}$
and $M$
allow to define the norm $\| df (x) \|$ of the differential $df$
at any point $x\in M\setminus Y$.

\begin{lem}
\label{lem-can-choose-unif-bounded}
One
can choose $f$ so that $\| df (x) \|$ is
uniformly bounded on $M\setminus Y$.
\end{lem}
We omit the proof of this technical result.

Fix a unitary trivialization of
$T M$ over $M\setminus Y$ as in the lemma above.
Pick an element $\widetilde{\phi}\in \tSymp_0\,
(M,\omega)$ and
represent it by a path $\{ \phi_t \}_{0\leq t\leq 1}$, $\phi_0 =
Id$, in $\Symp_0\, (M,\omega)$.
Consider the sets
\[
Y_{\{\phi_t\}} := \bigcup_{k\in\Z}\ \bigcup_{0\leq t \leq 1} \phi_t^{-k}
(Y),\ \ \
X_{\{\phi_t\}} := M\setminus Y_{\{\phi_t\}}.
\]
Since
$Y$ is of codimension at least 3, $Y_{\{\phi_t\}}$ has measure
zero
and therefore $X_{\{\phi_t\}}$ is a subset of full measure in $M$.

As it was stated in Remark~\ref{rem-lem-fuchs},
our trivialization of $TM$ over
$M\setminus Y$
symplectically identifies each tangent space of
$M\setminus Y$ with the standard symplectic space $\R^{2n}$.
Therefore one can view the differentials of the symplectomorphisms
$\phi_t$, $0\leq t\leq 1$, at any point of $X_{\{\phi_t\}}$ as
elements of the group ${\hbox{\it Sp}}\, (2n, \R)$. Thus to
any point $x\in X_{\{\phi_t\}}$ one can associate a path $\{
d\phi_t (x)\}_{0\leq t\leq 1}$ in ${\hbox{\it Sp}}\, (2n, \R)$
starting at the identity.
Denote by $F_{\{\phi_t\}}\,  (x)$ the element in
$\widetilde{\hbox{\it Sp}}\, (2n, \R)$ represented by this
path. Set $F_{\{\phi_t\}} = Id$ over $M\setminus
X_{\{\phi_t\}} = Y_{\{ \phi_t\} }$.

%\smallskip
\begin{lem}
\label{lem-integrability} For any path $\{\phi_t\}$ representing
$\widetilde{\phi}\in \tSymp_0\, (M,\omega)$
the function $\tau (F_{\{\phi_t\}} (\cdot) ) $ is integrable on
$M$.

\end{lem}
%\smallskip

\noindent {\bf Proof of Lemma~\ref{lem-integrability}.}

\noindent
For a point $x\in X_{\{\phi_t\} }$ denote
\[
V(x) := \widetilde{det}^2
(F_{\{ \phi_t\} } (x)).
\]
According to Proposition~\ref{prop-tau-barge-ghys}, $\tau$ is
continuous. Therefore in order to show that the function $\tau
(F_{\{\phi_t\} } (\cdot ))$ is integrable on $M$ it suffices to
show that $F_{\{\phi_t\} }: M\to \widetilde{Sp}\, (2n,
\R)$ has a bounded image in $\widetilde{Sp}\, (2n, \R)$.
In view of Proposition~\ref{prop-boundedness-sp-lagr-grassm} this
would follow if we show that the function
$V(x) := \widetilde{det}^2
(F_{\{ \phi_t\} } (x))$ is bounded on $X_{\{\phi_t\} }$.

The function $V$ can be described as follows.
Let $\Lambda\to M$ be the fiber bundle whose fiber over
a point $x\in M$ is the Lagrange Grassmannian $\Lambda
(T_x M)$ of the symplectic tangent space $T_x M$. Consider the
section $L : M\setminus Y\to \Lambda$ of the bundle $\Lambda\to M$
that associates to each $x\in M$ the Lagrangian plane
$L_0 \subset \Lambda (T_x M)= \Lambda(\R^{2n})$, where
$T_x M$ is identified with the standard symplectic linear
space $\R^{2n}$ by means of the symplectic (unitary)
trivialization of $T M$ over $M\setminus Y$ as above.
For a given
$x\in X_{\{\phi_t\}}$ the path $\{ \phi_t \}$ acting on $L (x)$
determines, by means of the trivialization, a path in $\Lambda
(\R^{2n})$ to which the map $\widetilde{det}^2$ associates a
number. This number is exactly $V (x)$.

Assume that $V (x)$ is not bounded. Then, since $F_{\{\phi_t\} }$
is continuous on $X_{\{\phi_t\}}$, there must exist a sequence $\{
x_k\}$ in $X_{\{\phi_t\}}$ which converges to a point $y\in
M\setminus X_{\{\phi_t\}}$ such that $\lim_{k\to +\infty} | V
(x_k) | = +\infty$. As $k\to
+\infty$, the paths $\gamma_k (t) := \phi_t
(x_k)$, $0\leq t\leq 1$, converge uniformly to a path
$\gamma_{lim} (t) := \phi_t (y)$, $0\leq t\leq 1$, lying in
$Y_{\{ \phi_t\} }$.

The total space $\Lambda$ of the bundle
$\Lambda\to M$ is compact. Therefore, possibly choosing a subsequence of $\{
x_k\}$, we may assume without loss of generality that the
Lagrangian subspaces $L (x_k)\in \Lambda (T_{x_k} M)\subset \Lambda$
converge to a Lagrangian subspace $L (y)\in \Lambda (T_y M)$ and
thus the paths $\{ d\phi_t (L (x_k) )\}_{0\leq t\leq 1}$ in
$\Lambda$ converge uniformly to a path
$\{ d\phi_t (L (y) )\}_{0\leq t\leq 1}$.

The total space ${\cal F}$ of the bundle ${\cal F}\to M$
is compact and the section
$f: M\setminus Y \to {\cal F}$,
defining our fixed unitary trivialization, has uniformly bounded
derivatives. Therefore, using Arzela-Ascoli theorem and possibly passing
to a subsequence, we may assume without loss of generality that
the sequence of maps $f\circ \gamma_k: [0,1]\to {\cal F}$
 converges $C^0{\hbox{\rm -uniformly}}$ to a {\it continuous} map
$g\circ \gamma_{lim}: [0,1]\to {\cal F}$, where $g$ is
a {\it continuous}
section of ${\cal F}\to M$ over the smooth
embedded curve
$\gamma_{lim} ([0,1])$.

Now each $V
(x_k)$ is the rotation number of the path
$\{ d\phi_t (L (x_k) )\}_{0\leq t\leq 1}$ of Lagrangian
planes in the bundle $\gamma_k^\ast TM$ over $[0,1]$ with respect
to the moving frame $f\circ\gamma_k (t)$ defining the
trivialization of that bundle.
Since the sequence $\{ L (x_k)\}$ converges to $L (y)$ and $\{ f\circ
\gamma_k\}$ converges uniformly to $g\circ \gamma_{lim}$, the
rotation numbers $V (x_k)$ converge to
a {\it finite} number $V_{lim} (y)$,
which is the rotation number of the path $\{ d\phi_t (L
(y) )\}_{0\leq t\leq 1}$ of Lagrangian planes
in the bundle $\gamma_{lim}^\ast TM$ over $[0,1]$ with respect to
the continuously moving frame $g\circ \gamma_{lim} (t)$ defining a
continuous trivialization of the bundle.
This
is in contradiction with our original assumption that $\lim_{k\to
+\infty} | V (x_k) | = +\infty$. The lemma is proven.\ \b
%\bigskip
%\bigskip

%\smallskip
\begin{lem}
\label{lem-corr-defined} Suppose that paths $\{\phi_t\}_{0\leq
t\leq 1}$, $\{ \psi_t \}_{0\leq t\leq 1}$ ,
$\phi_0 = \psi_0 = Id$,  in $\Symp_0\, (M,\omega)$ represent the
same element in $\tSymp_0\, (M,\omega)$. Then there exists
a subset $X_{\{\phi_t\}, \{\psi_t\} } \subset
X_{\{\phi_t\}}\cap X_{\{\psi_t \}}$ of the full measure in $M$ such that
for any $x\in X_{\{\phi_t\}, \{\psi_t\} }$
\[
F_{ \{\phi_t\} } (x) = F_{ \{\psi_t\} } (x) \in
\widetilde{\hbox{\it Sp}}\, (2n, \R).
\]
\end{lem}
%\smallskip

%\bigskip
\noindent {\bf Proof of Lemma~\ref{lem-corr-defined}.}

\noindent Since the paths $\{\phi_t\}$, $\{ \psi_t \}$ define the
same element in $\tSymp_0\, (M,\omega)$ there
exists a homotopy $\{ \varphi_{s,t}\}$, $0\leq s\leq 1$, $0\leq
t\leq 1$, between them which keeps the endpoints fixed. Consider
now the set
\[
Z_{\{\phi_t\}, \{\psi_t\} } := \bigcup_{k\in\Z}\ \bigcup_{0\leq s,t \leq 1}
\varphi_{s, t}^{-k} (Y).
\]
Since, according to its definition (see Lemma~\ref{lem-fuchs}),
$Y$ is a subset of codimension at least 3 in $M$, the set
$Z_{\{\phi_t\}, \{\psi_t\} }$ is of codimension at least 1 and therefore has measure
zero. The set
$X_{\{\phi_t\}, \{\psi_t\} } := M\setminus Z_{\{\phi_t\}, \{\psi_t\}
}$
is the one we have been looking for.
The lemma is proven.\ \b

%\smallskip
\begin{lem}
\label{lem-product-of-paths} Let $\{ \phi_t \}$, $\{ \psi_t \}$,
$0\leq t\leq 1$,
be paths in $\Symp_0\, (M,\omega)$ starting at the identity. Then there
exists a set $X^\prime_{\{ \phi_t\} , \{ \psi_t\} }$ of full measure in $M$
such that the following quantities are well-defined and
equal for any $x\in X^\prime_{\{ \phi_t\} , \{ \psi_t\} }$:
\begin{equation}
\label{eqn-F-phi-psi} F_{\{\phi_t\psi_t \}} (x) = F_{\{\phi_t\}
}\,  (\psi_1 (x))  F_{\{\psi_t \}} (x).
\end{equation}
\end{lem}
%\smallskip

%\bigskip
\noindent {\bf Proof of Lemma~\ref{lem-product-of-paths}.}

\noindent
Consider the set
\[
Z^\prime_{\{ \phi_t\} , \{ \psi_t\} } :=
\bigcup_{k\in\Z}\ \bigcup_{0\leq s,t \leq 1}
(\phi_t\psi_s )^{-k} (Y).
\]
Since $Y$ is of codimension at least 3 in $M$, the set
$Z^\prime_{\{ \phi_t\} , \{ \psi_t\} }$ has measure zero and
therefore
$X^\prime_{\{ \phi_t\} , \{ \psi_t\} } :=
M\setminus Z^\prime_{\{ \phi_t\} , \{ \psi_t\} }$
is a set of full measure in $M$.
Also
\[
X^\prime_{\{ \phi_t\} , \{ \psi_t\} } \subset X_{\{ \phi_t\psi_t\} }\cap
X_{\{ \phi_t\} } \cap X_{\{ \psi_t\} }
\]
and therefore for any
$x\in X^\prime_{\{ \phi_t\} , \{ \psi_t\} }$
both sides of the equality (\ref{eqn-F-phi-psi})
are well-defined.
For any such $x$
the element $F_{\{\phi_t\psi_t \}} (x) \in
\widetilde{\hbox{\it Sp}}\, (2n, \R)$ can be represented by a
path
\[
\{ d (\phi_t \psi_t) (x)\} = \{ d\phi_t (\psi_t (x)) \cdot d\psi_t
(x)\},
\]
$0\leq t\leq 1$, in ${\hbox{\it Sp}}\, (2n, \R)$.
But the
path $\{ d\phi_t (\psi_t (x))\}_{0\leq t\leq 1}$
is homotopic
with the fixed endpoints to the path $\{ d\phi_t
(\psi_1 (x))\}_{0\leq t\leq 1}$ in ${\hbox{\it Sp}}\, (2n, \R)$:
the homotopy is given by the 2-parametric family
\[
d\phi_t (\psi_s (x)),\ \ 0\leq t\leq 1, \ \ t\leq s\leq 1,
\]
of elements of
${\hbox{\it Sp}}\, (2n, \R)$,
which is well-defined since $\phi_t\psi_s (x) \notin Y$
for all $s,t$ and
$x\in X^\prime_{\{ \phi_t\} , \{ \psi_t\} }$.
Thus for any
$x\in X^\prime_{\{ \phi_t\} , \{ \psi_t\} }$ the paths
$\{ d\phi_t (\psi_t (x))\}_{0\leq t\leq 1}$ and
$\{ d\phi_t
(\psi_1 (x))\}_{0\leq t\leq 1}$ represent the same element in
$\widetilde{\hbox{\it Sp}}\, (2n, \R)$
and hence
\[
F_{\{\phi_t\psi_t \}} (x) = F_{\{\phi_t\} }\,  (\psi_1 (x))
F_{\{\psi_t \}} (x).
\]
The lemma is proven.\ \b

\subsection{Final steps in the construction of  ${\mathfrak f}$}

From this stage on the construction of ${\mathfrak f}$ proceeds in
exactly the same way as in \cite{Ba-Ghys}. Namely, let
$\widetilde{\phi}\in \tSymp_0\, (M,\omega)$
be represented by a path $\{ \phi_t \}_{0\leq t\leq 1}$, $\phi_0 = Id$,
in
$\Symp_0\, (M,\omega)$. For $x\in X_{\{\phi_t\}}$ and
$k\geq 1$ define
\[
\tau_k  (\{ \phi_t \}, x) = \tau \bigg( \prod_{i=0}^{k-1}
F_{\{\phi_t\} }\, (\phi_1^i (x)) \bigg).
\]
According to the definition of $X_{\{\phi_t\}}$, the number
$\tau_k  (\{ \phi_t \}, x)$ is
well-defined.
Since $\tau$ is a quasimorphism there exists a constant $C$ such
that for any $l$ and $x$
\begin{eqnarray}
\label{eqn-tau-k-tau-l} \tau_k (\{\phi_t \}, x) + \tau_l (\{\phi_t
\}, \phi_1^k (x)) - C \leq \tau_{k+l} (\{\phi_t \}, x) \leq
\nonumber\\ \leq \tau_k (\{\phi_t \}, x) + \tau_l (\{\phi_t \},
\phi_1^k (x)) + C.
\end{eqnarray}

%\smallskip
\begin{lem}
\label{lem-integr-tau-k} For any $k$ the function $\tau_k
(\{\phi_t \}, \cdot)$ is integrable on $M$.
\end{lem}
%\smallskip

%\bigskip
\noindent {\bf Proof of Lemma~\ref{lem-integr-tau-k}.}

\noindent Proceed by induction. For $k=1$ the statement follows
from Lemma~\ref{lem-integrability}. To prove the step of the
induction use (\ref{eqn-tau-k-tau-l}) with $l=1$ and
recall that $\phi_1$ is a symplectomorphism and therefore preserves
the measure.\ \b

%\smallskip
\begin{lem}
\label{lem-integral-tau-k-corr-defined} The integral $\int_M
\tau_k (\{\phi_t \}, \cdot)\,   d\mu$ does not depend on the choice
of the path $\{ \phi_t\}$ representing
$\widetilde{\phi}\in \tSymp_0\, (M,\omega)$.
\end{lem}
%\smallskip

%\bigskip
\noindent {\bf Proof of
Lemma~\ref{lem-integral-tau-k-corr-defined}.}

\noindent Let $\{ \phi_t\}$, $\{ \psi_t\}$ be two paths
representing the same $\widetilde{\phi}\in \tSymp_0\, (M,\omega)$.
Applying
Lemma~\ref{lem-corr-defined}
one gets that the functions $\tau_k (\{\phi_t \}, \cdot)$ and
$\tau_k (\{\psi_t \}, \cdot)$ differ on a set of measure zero in $M$
and
therefore have the same integral.
\b

\bigskip
\noindent Thus we can define $T_k (\widetilde{\phi}) := \int_M
\tau_k (\{ \phi_t \}, \cdot)\,   d\mu$ for any $\{ \phi_t\}$
representing $\widetilde{\phi}$. In view of
(\ref{eqn-tau-k-tau-l}),
\begin{equation}
\label{eqn-T-k-T-l} T_k (\widetilde{\phi}) + T_l
(\widetilde{\phi}) - C\mu (M) \leq T_{k+l} (\widetilde{\phi}) \leq
T_k (\widetilde{\phi}) + T_l (\widetilde{\phi}) + C\mu (M)
\end{equation}
for any $l$. Therefore  there exists a limit
\[
{\mathfrak f} (\widetilde{\phi}) := \lim_{k\to +\infty} T_k
(\widetilde{\phi}) /k.
\]
This is the definition of ${\mathfrak f}$. Now we will prove that
${\mathfrak f}$ is a non-trivial homogeneous quasimorphism.

%\bigskip
\subsection{
${\mathfrak f}$ is a quasimorphism}

Formula (\ref{eqn-T-k-T-l}) implies that
\begin{equation}
\label{eqn-T-T-1} |\ {\mathfrak f} (\widetilde{\phi}) - T_1
(\widetilde{\phi})\ | \leq C \mu (M)
\end{equation}
for some constant $C$ that does not depend on
$\widetilde{\phi}$. Now let $\{
\phi_t\}$, $\{
\psi_t\}$, $0\leq t\leq 1$,
be paths in $\Symp_0\, (M,\omega)$ starting at the identity and
representing, respectively, $\widetilde{\phi},
\widetilde{\psi}\in \tSymp_0\, (M,\omega)$.
The equality (\ref{eqn-F-phi-psi}) holds on a set
$X^\prime_{\{\phi_t\}, \{\psi_t \}}$ of full measure and therefore
\[
T_1 (\widetilde{\phi} \widetilde{\psi} ) = \int_M \tau \bigg(
F_{\{\phi_t\psi_t \} } (x) \bigg) d\mu = \int_M \tau \bigg(
F_{\{\phi_t\} }\,  (\psi_1 (x)) F_{\{ \widetilde{\psi}\} }\,  (x) \bigg)
d\mu.
\]
Since $\psi_1$ is a symplectomorphism,
\[
\int_M \tau \bigg( F_{\{\phi_t\} }\,  (\psi_1 (x)) \bigg) d\mu =
\int_M \tau \bigg( F_{\{\phi_t\} }\,  (x) \bigg) d\mu,
\]
and since $\tau$ is a quasimorphism,
\[
\bigg|\ \int_M \tau \bigg( F_{\{\phi_t\} }\,  (\psi_1 (x))
F_{\{\psi_t \}}\,  (x) \bigg) d\mu -
\int_M \tau \bigg( F_{\{\psi_t \}}\,  (x) \bigg) d\mu -
\]
\[
- \int_M \tau \bigg(
F_{\{\phi_t\} }\,  (x)
 \bigg) d\mu\ \bigg|
\leq C_1 \mu (M)
\]
for some constant $C_1$ independent of $\widetilde{\phi},
\widetilde{\psi}$. Therefore
\[
| T_1 (\widetilde{\phi} \widetilde{\psi} ) -  T_1 (\widetilde{\psi}) - T_1
(\widetilde{\phi}) | \leq C_1 \mu (M)
\]
for any $\widetilde{\phi}, \widetilde{\psi}$. Thus $T_1$ is a
quasimorphism. In view of (\ref{eqn-T-T-1}),
${\mathfrak f}$ is a quasimorphism as well.

%\bigskip
\subsection{
${\mathfrak f}$ is homogeneous}

Represent each $\widetilde{\phi}^k$, $1\leq k\leq m$, by a path
$\{ \phi^k_t\}_{0\leq t\leq 1}$, $\phi_0 = Id$. Observe
that, according to (\ref{eqn-F-phi-psi}),
\[
F_{\{\phi^m_t\}}\,  (x) =
\prod_{i=0}^{m-1} F_{\{\phi_t\} }\,
(\phi_1^i (x))
\]
for any $x$ from a set of full measure in $M$. For any such $x$
\[
\tau_k (\{\phi^m_t\}, x) =
\tau \bigg(
\prod_{i=0}^{k-1} F_{\{\phi^m_t\}}\, (\phi^{mi}_1 (x)) \bigg)  =
\tau \bigg(
\prod_{i=0}^{k-1} \prod_{j=0}^{m-1} F_{\{\phi_t\}}\,
(\phi^{mi+j}_1 (x))\bigg) =
\]
\[
= \tau \bigg(
\prod_{l=0}^{mk-1} F_{\{\phi_t\} }\,
(\phi^l_1 (x)) \bigg) =
\tau_{mk} (\{\phi^m_t\} , x)
\]
and therefore
\[
\ T_k
(\widetilde{\phi}^m ) = T_{mk} (\widetilde{\phi}).
\]
Hence
\[
{\mathfrak f} (\widetilde{\phi}^m ) = \lim_{k\to +\infty} T_k
(\widetilde{\phi}^m)/k =
 \lim_{k\to +\infty} T_{mk}
(\widetilde{\phi})/k = \]
\[
= m \cdot\lim_{k\to +\infty} T_{mk}
(\widetilde{\phi})/mk = m {\mathfrak f} (\widetilde{\phi} )
\]
and therefore ${\mathfrak f}$ is homogeneous.

\subsection{
${\mathfrak f}$ does not vanish on $\widetilde{\varphi}_{H_B}\in
\tHam  (M,\omega)$}

In order to prove that ${\mathfrak f}$ does not vanish on
$\tHam  (M,\omega)$ consider
the same Hamiltonian symplectomorphism that J.Barge and E.Ghys
used
to show that their homogeneous quasimorphism on $\Symp^{c}\, (B^{2n})$ does not vanish (see \cite{Ba-Ghys}, p.
263).

Namely, consider a ball $B\subset M^{2n}\setminus Y$.
Assume without loss of generality that our
chosen trivialization of $TM$ over $M\setminus Y$ coincides on $B$
with the trivialization defined by the identification of $B$
with a ball in the standard symplectic $\R^{2n}$.
Consider the Hamiltonian symplectomorphism $\varphi_{H_B}$
generated by $H_B$ and its lift $\widetilde{\varphi}_{H_B}$ in
$\tHam  (M,\omega)$. Following the
construction in \cite{Ba-Ghys} one easily sees that value of our
quasimorphism ${\mathfrak f}$ on $\widetilde{\varphi}_{H_B}$ is
the same as the value of the Barge-Ghys quasimorphism on
$\varphi_{H_B}$. According to \cite{Ba-Ghys}, the latter value is
non-zero. Hence ${\mathfrak f}$ does not vanish on
$\tHam  (M,\omega)$.

This finishes the proof of the claim that ${\mathfrak f}$ is a
homogeneous quasimorphism on $\tSymp_0\,
(M,\omega)$ that does not vanish on $\widetilde{\varphi}_{H_B}\in
\tHam  (M,\omega)$.
Theorem~\ref{thm-comm-length-closed-mfds} is proven.
\b

\bibliographystyle{alpha}

\end{document}